\documentclass[11pt]{article}
\usepackage{amssymb,latexsym}

\newtheorem{theorem}{Theorem}[section]

\newtheorem{defi}{Definition}[section]
\newtheorem{lemma}{Lemma}[section]
\newtheorem{prop}{Proposition}[section]

\def\binom#1#2{{#1}\choose{#2}}

\def\slfrac#1#2{\hbox{\kern.1em %
 \raise.5ex\hbox{\the\scriptfont0 #1}\kern-.11em %
 /\kern-.15em\lower.25ex\hbox{\the\scriptfont0 #2}}}

\newcommand{\eqn}[1]{(\ref{#1})}

\newcommand{\eeq}{\end{equation}}
\newcommand{\beql}[1]{\begin{equation}\label{#1}}
\newcommand{\bsq}{{\vrule height .9ex width .8ex depth -.1ex }}

\newcommand{\ZZ}{{\mathbb Z}}
\newcommand{\RR}{{\mathbb R}}

\newcommand{\bb}{{\bf b}}

\newcommand{\bbS}{\Sigma}

\newcommand{\sC}{{\cal C}}

\newcommand{\sE}{{\cal E}}

\newcommand{\sI}{{\cal I}}
\newcommand{\sJ}{{\cal J}}

\newcommand{\sM}{{\cal M}}

\newcommand{\sS}{{\cal S}}

\makeatletter
\def\@sect#1#2#3#4#5#6[#7]#8{\ifnum #2>\c@secnumdepth
     \def\@svsec{}\else
     \refstepcounter{#1}\edef\@svsec{\csname the#1\endcsname.\hskip .75em }\fi
     \@tempskipa #5\relax
      \ifdim \@tempskipa>\z@
        \begingroup #6\relax
          \@hangfrom{\hskip #3\relax\@svsec}{\interlinepenalty \@M #8\par}%
        \endgroup
       \csname #1mark\endcsname{#7}\addcontentsline
         {toc}{#1}{\ifnum #2>\c@secnumdepth \else
                      \protect\numberline{\csname the#1\endcsname}\fi
                    #7}\else
        \def\@svsechd{#6\hskip #3\@svsec #8\csname #1mark\endcsname
                      {#7}\addcontentsline
                           {toc}{#1}{\ifnum #2>\c@secnumdepth \else
                             \protect\numberline{\csname the#1\endcsname}\fi
                       #7}}\fi
     \@xsect{#5}}
\def\@begintheorem#1#2{\it \trivlist \item[\hskip \labelsep{\bf #1\ #2.}]}

\def\plain{plain}\ifx\fmtname\plain\csname fi\endcsname
     
     \input docstrip
     \preamble

     Do not distribute the stripped version of this file.
     The checksum in the header refers to the documented version.

     \endpreamble
     \generateFile{here.sty}{t}{\from{here.doc}{}}
     \endinput
\fi
\ifcat a\noexpand @\let\next\relax\else\def\next{%
    \documentstyle[here,doc]{article}\MakePercentIgnore}\fi\next
\ifx\@Hxfloat\@Hundef\else\expandafter\endinput\fi
\let\@Hxfloat\@xfloat
\def\@xfloat#1[{\@ifnextchar{H}{\@HHfloat{#1}[}{\@Hxfloat{#1}[}}
\def\@HHfloat#1[H]{%
\expandafter\let\csname end#1\endcsname\end@Hfloat
\vskip\intextsep\vbox\bgroup\def\@captype{#1}\parindent\z@
\ignorespaces}
\def\end@Hfloat{\egroup\vskip \intextsep}

\makeatother

\catcode`\@=11
\renewcommand{\section}{
        \setcounter{equation}{0}
        \@startsection {section}{1}{\z@}{-3.5ex plus -1ex minus
        -.2ex}{2.3ex plus .2ex}{\large\bf}
        }
\catcode`@=12

\begin{document}

\begin{center}
{\Large 
{\bf   Ternary Expansions  of Powers of $2$}
}\\

\vspace{1.5\baselineskip}
{\em Jeffrey C. Lagarias}\footnote{MSC Classification (2000):  11A63 (Primary), 11K16, 11K41,
26A18, 37A45 (Secondary)}\\
\vspace*{.2\baselineskip}
Dept. of Mathematics \\
University of Michigan \\
Ann Arbor, MI 48109-1109\\
\vspace*{1.5\baselineskip}

(To  Mel  Nathanson on his 60-th birthday)

\vspace*{2\baselineskip}
(July 11, 2008 ) \\
\vspace{3\baselineskip}
{\bf ABSTRACT}
\end{center}
P. Erd\H{o}s asked how frequently does $2^n$ have a   ternary expansion 
that omits the digit $2$.
He conjectured that this holds only for finitely many
values of $n$. We generalize this question to consider 
iterates of two discrete dynamical systems.
The first considers  truncated ternary expansions of real sequences 
$x_n(\lambda) =\lfloor \lambda 2^n \rfloor$,
where $\lambda  >0$ is a real number, along with
its untruncated version, while 
the second  considers  $3$-adic expansions
of  sequences 
$y_n(\lambda)= \lambda 2^n$, where $\lambda$ is
a $3$-adic integer. We show in both cases that the set of
initial values having
infinitely many iterates that omit the 
digit $2$ is  small in a suitable sense.  
For each nonzero initial value we obtain an 
asymptotic upper bound as $k \to \infty$  on the 
the number of the first $k$ iterates that omit the
digit $2$. We also study auxiliary problems concerning
the Hausdorff dimension of intersections of multiplicative
translates of $3$-adic  Cantor sets.

%
%
%
%
\setlength{\baselineskip}{1.0\baselineskip}

\section{Introduction}

P. Erd\H{o}s \cite{Er79} asked the question of when the
ternary expansion  of $2^n$ omits the digit $2$. This happens for $2^0 =(1)_3,$
$2^2=4 = (11)_3$ and $2^8=256= (100111)_3$. He conjectured that it
does not happen for all $n \ge 9$, and 
commented that: ``As far as I can see, there is no method at
our disposal to attack this conjecture.'' 
This question 
was initially studied by Gupta \cite{Gu78} who found by a sieving procedure that
there are no other solutions for $n < 4374$. In 1980 
Narkiewicz \cite{Na80} showed
 that the number 
 $$
 N_1(X) :=  \# \{n \le X: ~\mbox{the~ ternary ~expansion}~ (2^n)_3~\mbox{omits~the~digit}~2\}.
 $$
has $N_1(X) \le 1.62 X^{\alpha_0}$ with $\alpha_0 = \log_3 2 \approx 0. 63092$. 
The Erd\H{o}s  question remains open and has appeared in several problem lists,
e.g.  Erd\H{o}s and Graham \cite{EG80} and Guy \cite[Problem B33]{Guy2}.
In this paper we call the "Conjecture of Erd\H{o}s"   the weaker assertion that
there are only finitely many exponents $n$ such that the ternary expansion 
$(2^n)_3$ of $2^n$ omits the digit $2$.

This paper considers analogues of the conjecture of Erd\H{o}s for iterates of 
two discrete dynamical systems, one acting on the real numbers
and the other acting on the $3$-adic integers,  with an additional degree
of freedom given by a  parameter $\lambda$ specifying the initial condition. In both 
dynamical systems the parameter value  $\lambda=1$
recovers the original sequence $\{ 2^n: n \ge 0\}$ of Erd\H{o}s 
as a  forward  orbit of the dynamics.

The first dynamical system  is
$y \mapsto 2y$ acting on the real numbers, which is a
homeomorphism of $\RR$ that is an  expanding map.
It produces  a sequence of iterates $y_n= 2^n y_0$
starting from  $y_0= \lambda$. The {\em real dynamical system}
concerns the iterates $y_n$. We  also consider
an associated dynamical system which gives  integers, by  applying the
floor operator, obtaining
the sequence   $x_n= \lfloor y_n \rfloor$; that is, 
\beql{101a}
x_n = x_n(\lambda) :=  \lfloor \lambda 2^n \rfloor,~~~~~\mbox{for}~~ n \ge 0.
\eeq
We call this the {\em truncated real dynamical system}. Strictly speaking  the 
truncated real dynamical system has  forward orbits 
involving two variables $O^{+}(\lambda):= \{ (y_n(\lambda), x_n(\lambda)): n \ge 0\}$,
with $\{ y_n(\lambda)\}$ driving the dynamics. However
the expanding nature of the map $y \mapsto 2y$ implies that
the integer sequence  $\{x_n(\lambda): n \ge 0\}$ 
 contains enough information 
to uniquely determine the initial condition $\lambda$ of the
iteration; here we consider the ternary expansions of the $x_n(\lambda)$.

The second dynamical system is $y \mapsto 2y$ acting on the
$3$-adic integers $\ZZ_3$, which is a  3-adic measure-preserving 
homeomorphism of $\ZZ_3$.
It produces a sequence of iterates  $y_n = 2^n y_0$ 
starting from the initial condition 
$y_0= \lambda$. We write
\beql{101b}
y_n = y_n(\lambda) = \lambda 2^n,~~~~~\mbox{for}~~ n \ge 0,
\eeq
In this case we study membership of values $ y_n(\lambda)$
in the subset
$\Sigma_{3, \bar{2}}$ of all $3$-adic integers whose $3$-adic
expansion omits the digit $2$; this is the multiplicative translate
$\frac{1}{2} \Sigma_{3, \bar{1}}$ of the $3$-adic analogue
$\Sigma_{3, \bar{1}}$ of the classical "middle-third" Cantor set.

In the real number case dynamical systems of a
 related nature  have been  studied by
 several authors.
 Flatto, Lagarias and Pollington \cite{FLP95} introduced 
a parameter $\lambda$ in similar
questions concering the fractional parts  
of the sequences $\{\{ \lambda \xi^n\}\}$, for fixed $\xi>1$,
with the aim of proving results for the
parameter value $\lambda=1$ by proving universal results
valid for all parameter values $\lambda >0$. Recently 
Dubickas and Novickas \cite{DN05}
considered the prime or compositeness properties
of integers occurring in truncated recurrence  sequences, including
$\lfloor \lambda 2^n \rfloor$ as a particularly simple case.
Dubickas \cite{Du06} further extends both these results
to  certain $\lambda$ that are  real algebraic numbers. 

The paper contains both results and conjectures;.
We now state them in detail. 

%
%
%
%

%
%
%
%
\subsection{Truncated Real Dynamical System: Results}

For the truncated real dynamical system $x_n = \lfloor \lambda 2^n\rfloor$, 
we  show that  there is a uniform asymptotic upper bound valid for all nonzero
$\lambda$ on the number of
$n \le X$ for which $(\lfloor \lambda 2^n\rfloor)_3$ omits the digit $2$.
 Let $(k)_3$ denote the ternary digit expansion of the integer $k$.

%
%

\begin{theorem}~\label{th11}
For each $\lambda >0$, the upper bound 
\beql{102}
N_{\lambda}(X):= \# \{n : ~ 1\le n \le X~\mbox{and}~ 
(\lfloor \lambda 2^n\rfloor)_3~ \mbox{omits~the~digit~2} \} \le 25 X^{0.9725}
\eeq
holds  for all all sufficiently large $X \ge n_0(\lambda).$
\end{theorem}

In the complementary direction,  the function $N_{\lambda}(X)$
is not always bounded. The next result shows there  exist 
uncountably many $\lambda >0$ such that the sequence $x_n(\lambda)$
contains infinitely many integers omitting the digit $2$
in their ternary expansion. 
%
%
\begin{theorem}~\label{th12}
There exists an infinite sequence $S=\{n_k: k \ge 1\}$
satisfying $n_1=2$ and
\beql{103} 
2^{\frac{1}{14}(n_{k-1} +2k -7)} \le n_{k} \le 2^{27(n_{k-1} +2k + 6)},
\eeq
having the following property: The  set of real numbers $\Sigma(S)$
consisting of all $\lambda >0$
for which all the integers 
$\{ x_{n}(\lambda):=\lfloor \lambda 2^{n} \rfloor: n \in S\}$
have  ternary expansions omitting the digit $2$ is an uncountable set.
\end{theorem}

The set of exponents produced in this
theorem forms  a very thin infinite set.
One can show that \eqn{103} implies that for $X \ge 2$, its cardinality 
satisfies 
\beql{103b} 
   \# \{ n_k: ~1 \le n_k \le X \} \ge  \log_{\ast} (X )- 4.
\eeq
in which $\log_{\ast}(X)$
denotes the number of iterations of the logarithm function
starting at $X$ necessary to get a value of smaller than $1$.
Thus we obtain that for all $\lambda \in \Sigma(S)$, 
\beql{103c}
N_{\lambda}(X) \ge \log_{\ast} (X) - 4.
\eeq

We next consider properties of the set of $\lambda$ that have infinitely such integers. 
We define  the {\em truncated real exceptional set}  $\sE_{T}(\RR_{+})$ by
\beql{104} 
 \sE_{T}(\RR_{+}) := \{ \lambda > 0: ~ \mbox{infinitely~many~ternary~expansions}~
(\lfloor \lambda 2^n\rfloor)_3 
~\mbox{omit~the~digit}~2 \}
\eeq
We prove the following result.

%
%
\begin{theorem}~\label{th13}
The truncated real exceptional set 
has Hausdorff dimension
$$
\dim_{H}(\sE_{T}(\RR_{+}))= \log_3 (2 )= \frac{\log2}{\log3}\approx  0. 63092.
$$
It has nonzero $\log_3 (2)$-dimensional Hausdorff measure.
\end{theorem}

This result gives an indication why it may be a hard problem to tell
whether there are infinitely many exceptional powers of $2$
for any particular $\lambda$, such as $\lambda=1$. Namely, it is likely 
to be a hard problem  to  decide whether any particular real number
belongs  to this "small" exceptional set.\\

%
%
%
%
\subsection{Real Dynamical System: Conjecture }

Consider  the real dynamical system $y \mapsto 2y$ on $\RR_{+}.$
without truncation, having  forward orbits $O^{+}(\lambda):=\{ y_n = \lambda 2^n: n \ge 0\}$. 
We define the  {\em  real exceptional set} $\sE(\RR_{+})$ by
\beql{104bb}
\sE(\RR_{+}) := \{ \lambda > 0: ~ \mbox{infinitely~many~ternary~expansions}~
( \lambda 2^n)_3 
~\mbox{~omit~the~digit}~2 \}.
\eeq
This set is much more constrained than the truncated exceptional set $\sE_T(\RR_{+})$
discussed above.
As far as we know it could
even be the empty set. The conjecture of Erd\H{o}s  is equivalent to the assertion that
$1 \not\in \sE(\RR_{+}).$

Concerning  this exceptional set  we make the following conjecture.

%
%
\paragraph{Conjecture A.}
{\em The  real  exceptional set 
$$
\sE(\RR) := \{ \lambda \in \RR_{+}:~ \mbox{infinitely many ternary
expansions}~ (\lambda 2^n)_3~ \mbox{omit 
the digit 2}\}
$$
has Hausdorff dimension zero.} \\

A stronger form of this conjecture would be that the exceptional set is countable;
even stronger would be the assertion that the real exceptional
set is empty. Thus, for the moment, there remains the possibility that 
the conjecture of Erd\H{o}s might hold for all initial conditions $\lambda >0$,
for the full ternary expansions $(\lambda 2^n)_3$ as real numbers. 

Note that if  the real exceptional set is nonempty, it will necessarily be an infinite set, because it
is forward invariant under multiplication by $2$, i.e.
$2 \sE(\RR_{+}) \subset \sE(\RR_{+}).$ It is clearly also forward invariant
under multiplication by $3$, i.e.
$3 \sE(\RR_{+} )\subset \sE(\RR_{+})$. Thus it is forward invariant under two commuting
semigroup actions.  But the real exceptional set  is not known to be a (topologically) closed set,
so that results on Hausdorff dimension on closed sets invariant under commuting
semigroup actions cannot be directly applied. 


%
%
%
%
\subsection{$3$-Adic  Dynamical System: Results}

For a $3$-adic integer $\lambda= \sum_{j=0}^{\infty} d_j 3^j$ with each $d_j \in \{0,1,2\}$
we write $(\lambda)_3 = ( \cdots d_2 d_1d_0)_3$ for its $3$-adic digital expansion.
Our first observation is an upper bound on the number of solutions valid
for all nonzero $\lambda \in \ZZ_3$, which extends the result
of Narkiewicz \cite{Na80}
for  $\lambda=1$, using essentially  the same proof.

%
%
\begin{theorem}~\label{th14}
For each nonzero $\lambda \in \ZZ_3$, the $3$-adic integers,
and each $X \ge 2$, 
\beql{107}
\tilde{N}_{\lambda}(X):=\# \{ n \le X:~ (\lambda 2^n)_3 \in \ZZ_3 ~\mbox{omits~the~digit}~2 \} 
\le 2X^{\alpha_0},
\eeq
with $\alpha_0= \log_3 2 \approx 0.63092$.
\end{theorem}

We next study 
the {\em $3$-adic exceptional set}
\beql{109}
 \sE(\ZZ_3) := \{ \lambda \in \ZZ_3: ~ \mbox{infinitely~many~ 3-adic~expansions}~
\lambda 2^n  ~\mbox{omit~the~digit~2} \}.
\eeq
This set seems hard to study directly, so as
approximations to the $3$-adic exceptional set, 
we define for $k \ge 1$ the sequence of sets 
\beql{109b}
 \sE^{(k)}(\ZZ_3) := 
\{ \lambda \in \ZZ_3: ~ \mbox{at~least~ $k$~ values~ of}~
\lambda 2^n  ~\mbox{omit~the~digit~2} \}.
\eeq
These sets clearly form a nested family under inclusion, 
$$
\sE^{(1)}(\ZZ_3) \supset \sE^{(2)}(\ZZ_3) \supset \sE^{(3)} (\ZZ_3) \supset \cdots,
$$
and their intersection contains the exceptional set $\sE(\ZZ_3).$
These sets are somewhat easier to study.

We consider the problem of estimating the Hausdorff dimension
of the sets $ \sE^{(k)}(\ZZ_3)$ (with respect to the $3$-adic metric) and show the following result.

%
%
\begin{theorem}~\label{th15}
(1) The exceptional set $\sE^{(1)}(\ZZ_3))$ has Hausdorff dimension
\beql{170}
\dim_H( \sE^{(1)} (\ZZ_3)) = \alpha_0 \approx 0.63092.
\eeq

(2) The exceptional set $\sE^{(2)}(\ZZ_3)$ has Hausdorff dimension bounded by
\beql{171}
 \frac{1}{2} \log_3 (2) \le \dim_H (\sE^{(2)}(\ZZ_3)) \le  \frac{1}{2}.
\eeq

(3) The exceptional set $\sE^{(3)}(\ZZ_3)$ has positive Hausdorff dimension bounded by
\beql{172}
\frac{1}{6} \log_3 2 \le  \dim_H (\sE^{(3)}(\ZZ_3)) \le \dim_H (\sE^{(2)}(\ZZ_3)). 
\eeq
\end{theorem}

This result is only a beginning of the study of $dim_{H} (\sE^{(k)})$
for general $k$.
The (not necessarily closed) set 
 $\sE^{(k)}(\ZZ_3)$
 is a countable union of closed sets 
 $\sC(2^{m_1}, 2^{m_2}, \cdots, 2^{m_k})$ consisting of those $\lambda$ for which 
 $\{\lambda 2^{m_j}: 1 \le j \le k\}$  all have 3-adic expansions
 that omit the digit $2$. One can use this  to obtain upper and lower
 bounds on Hausdorff dimension of these sets by analyzing the Hausdorff dimension
 of the individual sets $\sC(2^{m_1}, 2^{m_2}, \cdots, 2^{m_k}).$
 These sets are 
 intersections of
multiplicative translates of the $3$-adic Cantor set, which
we discuss in the next subsection. 
In Theorem~\ref{th15}  the upper bound
in (2) is deduced using  Theorem~\ref{th16} below. 

It is not clear whether $\dim_H (\sE^{(k)}(\ZZ_3)) >0$ for all $k \ge 1$.
Proving or disproving this assertion already seems a subtle question.

Since $\sE(\ZZ_3) \subseteq \sE^{(k)}(\ZZ_3)$ for each $k \ge 1$,  
 any upper bound on
the Hausdorff dimension of $\sE^{(k)}(\ZZ_3)$ gives
an  upper bound
for the Hausdorff dimension of the 
 $3$-adic exceptional set $ \sE(\ZZ_3).$
Each condition $\lambda 2^{m_j}  \in \Sigma_{3, \bar{2}}$
imposes  more constraints, apparantly lowering the 
Hausdorff dimension. This motivates the following
conjecture concerning the $3$-adic exceptional set $\sE(\ZZ_3).$ 

%
%
\paragraph{Conjecture B.}
{\em The   $3$-adic exceptional set 
$$
\sE(\ZZ_3) := \{ \lambda \in \ZZ_3: \mbox{infinitely  many 3-adic
expansions} ~\lambda 2^n ~ 
\mbox{omit the digit 2}\}
$$
has Hausdorff dimension zero.} \\

As in the real dynamical system case, we do not know much about this exceptional set, except that
it contains $0$. Again, the conjecture of
Erd\H{o}s is equivalent to the assertion that $1 \not\in \sE(\ZZ_3)$.
The $3$-adic exceptional set $\sE(\ZZ_3)$
 is forward invariant under multiplication by $2$ and multiplication
by $3$, but is not known to be a closed set. \\

%
%
%
%

\subsection{Intersection of Multiplicative Translates of Cantor Sets: Results}

The study of the exceptional sets $\sE^{(k)}(\ZZ_3)$ leads to
auxiliary questions concerning the Hausdoff dimensions of
intersections of 
multiplicative translates of the standard $3$-adic Cantor set $\Sigma_{3, \bar{2}}$, defined by
\beql{160}
\Sigma_{3, \bar{2}} :=\{ \lambda \in \ZZ_3:~\mbox{the ~3-adic~expansion} ~ (\lambda)_3~
\mbox{omits~the~digit~2} \}.
\eeq
For integers $1 \le M_1< M_2 < \cdots < M_k$ we  study the {\em multiplicative intersection sets}
\begin{eqnarray}\label {161}
\sC(M_1, M_2, \cdots , M_k) & := & \{ \lambda \in \ZZ_3: ~ (M_j \lambda)_3 ~~\mbox{omits~the~digit}~2
\mbox{ for}~ 1 \le j \le k\} \nonumber \\
&=& \bigcup_{j=1}^k  \left(\frac{1}{M_j }\Sigma_{3, \bar{2}}\right)
\end{eqnarray}
These sets are  closed sets.
The standard "middle third" Cantor set 
\beql{160b}
\Sigma_{3, \bar{1}} := \{ \lambda \in \ZZ_3:~\mbox{the~3-adic~digit~expansion}~(\lambda)_3~
\mbox{omits~the~digit~1} \}.
\eeq
has $\Sigma_{3, \bar{1}}= 2 \Sigma_{3, \bar{2}}$,  so that all results given below for $\Sigma_{3 \bar{2}}$
convert to  equivalent results for multiplicative translates of $\Sigma_{3, \bar{1}}.$

Multiplicative intersection sets arise in studying sets $\sE^{(k)}(\ZZ_3)$,
because they are given by countable unions of such sets, namely
$$
\sE^{(k)}(\ZZ_3) = \bigcup_{ 0 \le m_1< m_2< ...< m_k}  \sC(2^{m_1} , 2^{m_2}, \cdots, 2^{m_k})
$$

What can be said about the Hausdorff dimension of sets
$\sC(M_1, M_2, ..., M_k)$? This dimension depends in
a complicated manner on the $3$-adic expansions of the $M_i$,
and leads to various  problems which seem interesting in their
own right. 

%

\begin{theorem}~\label{th16}
Let $M$ be a positive integer which is not a power of $3$. Let $\Sigma_{3, \bar{2}}$
be the ternary Cantor set. 
Then the Hausdorff dimension of
$\sC(1, M) = \Sigma_{3, \bar{2}} \cap \frac{1}{M} \Sigma_{3, \bar{2}}$ satisfies
\beql{181}
\dim_{H} (\sC(1, M)) \le \frac{1}{2}.
\eeq
\end{theorem}

 We do not know if this bound is sharp.  However it is possible to show that 
 $$
 \dim_{H} (\sC(1, 7)) = \log_3 ( \frac{1+ \sqrt{5}}{2}) \approx 0.438.
 $$

For lower bounds on the Hausdorff dimension of such sets, we give the following
sufficient condition for positivity of the  Hausdorff dimension. 
%

\begin{theorem}~\label{th17}
Let $1 \le M_1< M_2 < \cdots < M_k$ be positive integers. Suppose there is
a positive integer $N$ belonging to the $3$-adic Cantor set
$\Sigma_{3, \bar{2}}\cup \ZZ$ such that  all the integers $NM_i$ satisfy
\beql{182}
N M_i \in \Sigma_{3, \bar{2}} \cap \ZZ,~~ 1 \le j \le k.
\eeq
Then 
\beql {183}
dim_{H}(\sC(M_1, M_2, ..., M_k) )\ge  \frac{\log_3 (2)}{\lceil \log_3 (NM_k)\rceil }.
\eeq
\end{theorem}

This is proved by direct construction of a Cantor set 
of positive Hausdorff dimension inside $\sC (M_1, M_2, ..., M_k)$.  

This result  gives a possible approach to obtaining a nonzero lower bound
for $\dim_{H}(\sE^{(k)}(\ZZ_3))$ for $k=4$ or larger, if suitable $M_i=2^{n_i}$ can
be found that fulfill its hypotheses.  However it can be shown that the sufficient
condition of Theorem \ref{th17} is not necessary, e.g.  $N=1$ and 
$M_1=1, M_2=52$ does not satisfy the hypothesis of this theorem, but
$\sC(1, 52)$ has positive
Hausdorff dimension. Thus  further strengthenings of this approach may be  possible.


Determining the structure and Hausdorff dimension of the sets $\sC(M_1, ..., M_k)$
leads to many open problems.\\

{\em Problem 1.} Let 
$$
\sM_{C} := \{ M \ge 1: ~\mbox{there~exist~integers} N_1, N_2 \in \Sigma_{3, \bar{2}}~
\mbox{with} ~N_1M= N_2\}.
$$
Obtain upper and lower bounds for the number
of integers $1 \le M \le X$ in $\sM_{C}$.\\

{\em Problem 2.} Let 
$$
\sM_{H}:= \{ M \ge 1: ~\dim_{H}( \sC(1,M) >0.\}
$$
Obtain upper and lower bounds for the number
of integers $1 \le M \le X$ in $\sM_{H}$.\\

These are different problems, because
it  can be shown that the inclusion  $\sM_{C} \subset \sM_{H}$  is strict.

%
%
%
%

\subsection{Generalization of  the Erd\H{o}s  Conjecture }

We formulate the   following strengthening of Erd\H{o}s's original question,
by analogy with  a conjecture of Furtstenberg \cite[Conjecture 2']{Fu70},
which is reviewed  in  \S5.

%
%
\paragraph{Conjecture E.} {\em Let $p$ and $q$ be multiplicatively independent
positive integers, i.e. all $\{ p^iq^j: i \ge 0, j\ge 0\}$ are distinct. Then
the base $q$ expansions of the powers $\{ (p^n)_{q}: n \ge 1\}$ have the property that 
any given finite pattern  $P=a_1 a_2\cdots a_k$ of consecutive $q$-ary digits occurs in
$(p^n)_{q}$,  for all sufficiently large  $n \ge n_0(P)$.} \\

Conjecture E generalizes Erd\H{o}s's original problem, which  is the special case $p=2$, $q=3$
with the single  pattern $P=2$.
We note that  Furstenberg's original conjecture concerns $d$-ary expansions 
of $\{ (p^n)_{d}: n \ge 1\}$ with $d=pq$ 
in which $p$ and $q$ are multiplicatively independent, i.e. his conjecture would apply to 
the $6$-adic expansion  $\{ (2^n)_6: n \ge0 \}$, rather than the $3$-adic 
expansion above.

This conjecture  might more properly be formulated  as
a question, since we present no significant new evidence in its favor.
However we think that any mechanism that forces a single pattern to appear 
from some point on should apply to all patterns.

%
%
%
%

\subsection{Summary }

First, this paper
places the original Erd\H{o}s problem in a more general dynamical context. 

The  two dynamical generalizations  seem to give restrictions on the original Erd\H{o}s question
 of  roughly equal strength, as formulated in Theorems \ref{th11} and \ref{th14}.
 That is, they each reduce the number of candidate $1 \le n \le X$ to at most
 $X^{c}$ for some $0< c<1.$
What is interesting is that these arguments  use  "independent"
information about the ternary expansions of $2^n$. 
The method used 
for the real dynamical system 
estimates the omission
of $2$ in the $\log_3 X$ most significant
ternary digits of $2^n$, while for the $3$-adic dynamical system 
the method estimates the 
omission of $2$ in the $\log_3 X$ least 
significant ternary digits of $2^n$. Heuristically, the most
significant digits and least significant digits seem uncorrelated;
this is the "independence" referred to above.
Furthermore, since the ternary expansion 
$(2^n)_3$ has about $\alpha_0 n$
ternary digits, the vast number of digits in 
the middle of the expansion are not exploited in either method;
only a logarithmically
small proportion of the available digits in the 
ternary expansion   $(2^n)_3$ are considered  in the two methods. 

It seems  a challenging problem to find a method that effectively
combines the two approaches to find better upper bounds on
$N_1(X)$ than that given by Narkiewicz. Can one obtain an upper
bound of $O( X^{\beta})$ for some $\beta < \log_3 2$ in this way?
Can one show that the high order digits and the low order
digits in the ternary expansion $(2^n)_3$ are
"uncorrelated" in some quantifiable way? 

Second, we  formulate
 Conjecture A and Conjecture B , asserting  Hausdorff dimension zero of 
 exceptional sets, which  seem
 more approachable questions than
the original question of Erd\H{o}s. 
A much harder question seems to be to resolve whether the 
exceptional sets $\sE(\RR_{+})$ and $\sE(\ZZ_3)$ are countable or finite. 

Third, our analysis  leads to a variety of interesting auxiliary problems in
 combinatorial number theory. These  
 concern the Hausdorff dimension  of intersections of multiplicative
translates of $3$-adic Cantor sets. These Hausdorff dimensions depend in
an complicated arithmetic way on the values of the integer multipliers. 
These sets seem worthy of further study. 

Finally, we observe analogies with work of Furstenberg \cite{Fu67}, \cite{Fu70}
on actions of multiplicative semigroups and intersections of Cantor sets.
This resulted in formulating  Conjecture E.

%
%
%
%

\subsection{Contents and Notation}
 
The contents of the rest of the paper are as follows.
In \S2 we prove results
for the truncated real dynamical system. In \S3 we prove
results for the $3$-adic dynamical system.
In \S4 we establish  auxiliary results on  the Hausdorff dimensions of 
intersections of a finite number of multiplicative translates  (by positive integers)
of the $3$-adic  Cantor set, and include several examples. 
 These results are used to complete the
proofs of one result in \S3. In \S5
we discuss work of Furstenberg. This includes
a conjecure  which motivates Conjecture E, 
and his formuation of a notion  transversality of semigroup actions
on a compact space and  implications for intersections of Cantor sets. 
In the concluding  section \S6  we describe history 
associated to Erd\H{o}s's  original question.


\paragraph{Notation.} Let
$$\{\{x\}\} := x - \lfloor x \rfloor = x~(\bmod~ 1)$$
denote the fractional part of a real number $x$.
Let
$$
\langle \langle x \rangle \rangle := \{\{x+1/2\}\} - 1/2
$$
denote the (signed) distance of $x$ to the nearest integer.

\paragraph{Acknowledgments.}
I am grateful to  A. Pollington, 
K. Soundararajan and H. Furstenberg for helpful
comments and references.  I thank the reviewer for helpful comments
and suggestions. 
The author was supported
by NSF grant DMS-0500555.

%
%
%
%
\section{Real Dynamical System: Proofs}

We consider the sequence of real numbers $x_n^{\ast}:= \lambda 2^n$,
and consider the associated integers 
$$x_n(\lambda) = \lfloor x_n^{\ast} \rfloor.$$
On taking logarithms to base $3$ we have
$$ 
\log_3 x_n^{\ast} = \log_3 \lambda + n \log_3 2 = m_n + w_n,
$$
in which $m_n = \lfloor \log_3 x_n^{\ast} \rfloor$ is the integer part and 
$ w_n :=\log_3 x_n^{\ast} ~(\bmod~1)$ is the fractional part, with 
$0 \le w_n < 1$.  Now  the digits in the 
ternary expansion of  $ x_n(\lambda)$  are completely determined by knowledge of 
the real number $w_n$, since $x_n(\lambda)= 3^{m_n}  3^{w_n}$, 
so they are the first $m_n$ ternary digits in the ternary expansion of $3^{w_n}$, 
since multiplication by $3^{m_n}$ simply shifts ternary digits to the left
without changing them.

On the other hand, the sequence of $w_n$ form an orbit under
iteration of  the map 
$T: [0,1] \mapsto [0,1]$ given by
\beql{201}
T(w) = w +\log_3 2 ~(\bmod ~1).
\eeq
on taking initial condition $w_0=  \log_3 \lambda$, with 
  $w_{n+1}= T(w_n)$.
Since $\alpha_0=\log_3 2$ is irrational, the map $T$ is an irrational
rotation on the torus $\RR/\ZZ$, which is known to be uniquely ergodic. In particular, 
every forward orbit of iteration of $T$ is uniformly distributed
$(\bmod~1)$, with the convergence rate to uniform distribution
determined by properties of the continued fraction expansion of $\alpha_0$.
We now examine the consequences of this property  for the ternary
expansions of $x_n^{\ast}$.

First, the leading ternary digits of $3^{w_n}$ specify  the position 
of $w_n$ in the interval $[0,1]$ to a small subinterval.  The property of omitting
the digit $2$ in a leading digit of a ternary expansion of $x_n$ will
prohibit $w_n$ from  certain subintervals  in $[0,1];$ the allowed subintervals
will have small measure. 
Using the fact that the distribution of $w_n (\bmod~1)$ approaches
the uniform distribution  fairly rapidly, one  can show that most $w_n$
have some leading digit that is a $2$; Theorem~\ref{th11} is deduced
using this idea, where the number $k$ of leading digits used will 
depend on the interval $[1, X]$ considered. 

Second, one use a construction selecting
a rapidly growing set of values of $n=n_k$, chosen using the
continued fraction expansion of $\alpha_0$,  in  such a way as to
permit each  $w_{n_k}$ to fall in a "good" interval where the initial ternary
digits for a large set of short  intervals 
have $x_{n_k}(\lambda)$'s 
with ternary expansions avoiding any $2$'s. 
A recursive intervals construction, which modifies
$\lambda$ slightly at each stage while not disturbing the initial ternary
digits already selected,  produces
the sets in Theorem~\ref{th12}. Finally, we use  a quantitative version
of such an intervals construction producing the set of Hausdorff dimension $\alpha_0$
in Theorem~\ref{th13}. 

We begin with two preliminary lemmas, the first on the spacings of
multiples of an irrational number (modulo one) and the second on
Diophantine
approximation properties of $\alpha_0= \log_3 2$.

%
%

\begin{lemma}~\label{lem21}
Let $\theta$ be irrational and consider the $N+1$ numbers 
$$\{ x+ j \theta ~(\bmod ~1): 0 \le j \le N\},$$
viewed as subdividing the 
torus $\RR/ \ZZ$  (the interval $[0,1]$ with endpoints
identified) into $N+1$ subintervals ("steps").

(1) These subintervals 
take at most three distinct lengths. If three
different lengths occur, say $L_1, L_2, L_3$,  then one of them is
the sum of the other two, say $L_1+L_2=L_3$.

(2) Let the continued fraction expansion of $\theta=[a_0, a_1, a_2, \cdots]$,
have partial quotients $a_i$ and convergents  $\frac{p_n}{q_n}$
with denominators  satisfying
$q_{n+1} = a_{n+1} q_n + q_{n-1}.$ Write uniquely 
\beql{210a}
N= (j+1) q_n + q_{n-1} + k,~~~ 0 \le k \le q_n -1
\eeq
with $0 \le j \le a_{n+1} -1.$ Then the subintervals have lengths
\begin{eqnarray*}
L_1  &=  & | \langle \langle q_n \theta \rangle \rangle|  \\
L_2  &=&   |  \langle \langle q_{n-1}\theta \rangle \rangle  + 
(j+1) \langle \langle q_n \theta \rangle \rangle|\\
L_3  &=&   |  \langle \langle q_{n-1}\theta \rangle \rangle  +  j \langle \langle q_n \theta \rangle  \rangle|
\end{eqnarray*}
and occur with  multiplicities $jq_n+q_{n-1}+k+1, ~k+1,$ and $~ q_n-(k+1),$ respectively. 
Here $L_3=L_1+L_2$, and $L_1< L_2$ if $0 \le j \le a_{n+1}-2$, while $L_2< L_1$
if $j=a_{n+1}-1$.  The intervals of  size $L_3$ do not occur if and only if $k= q_{n}-1$.

(3) For $N= q_{n+1}-1$, there occur intervals of exactly two lengths $L_1, L_2$
as above, and these lengths satisfy
\beql{212a}
L_2 < L_1 < 2L_2.
\eeq
\end{lemma}

\paragraph{Proof.} 
(1), (2)  These results have a long history, which is detailed in Slater \cite{Sl67}.
In particular, (2)  implies (1) and the formulas in (2) appear in Slater \cite[eqn. (33), p. 1120]{Sl67}.
The ordering of $L_1$ and $L_2$ follows from the fact that the $\langle\langle q_n \theta\rangle \rangle$
alternate in sign with successive $n$.

(3) Let $N=q_{n+1}-1$. If $a_n \ge 2$ then the decomposition \eqn{210a} is 
$$
N= (a_{n+1}-1) q_n + q_{n-1} + (q_n-1)
$$
 with $k= q_n-1$ and $j=a_{n+1}-1$,
Now (2) says  there are steps of exactly 
two lengths $L_1$ and $L_2$ given by
\begin{eqnarray*}
L_1  &=  & | \langle \langle q_n \theta \rangle \rangle|  \\
L_2  &=&   |  \langle \langle q_{n-1}\theta \rangle \rangle  + 
(a_{n+1}-1) \langle \langle q_n \theta \rangle \rangle|
\end{eqnarray*}
and $L_2< L_1$. Next we have 
$$
  \langle \langle q_{n+1} \theta \rangle \rangle=
\langle \langle q_{n-1}\theta \rangle \rangle  + 
a_{n+1} \langle \langle q_n \theta \rangle \rangle =
(\langle \langle q_{n-1}\theta \rangle \rangle  + 
(a_{n+1}-1) \langle \langle q_n \theta \rangle \rangle) +( \langle \langle q_n \theta \rangle \rangle).
$$
Since $\langle \langle q_{n+1}\theta \rangle \rangle$ and 
$ \langle \langle q_{n}\theta \rangle \rangle$ have opposite signs, and
$$
|\langle \langle q_{n+1}\theta \rangle \rangle| \le L_2 
$$
we must have
$$
L_2 < L_1 = L_2 + |\langle \langle q_{n+1}\theta \rangle \rangle| < 2L_2.
$$
(The fact that $\theta$ is irrational gives the strict inequality at the last step.) 

There remains the case $a_{n+1}=1$.
Now we find that the decompostion \eqn{210a} is 
$$
N= q_{n}+ q_{n-1}-1 = a_n q_{n-1} + q_{n-2} + (q_{n-1}-1),
$$
 with $k= q_{n-1}-1$ and $j=a_{n-1}-1$.
As before, there are intervals of exactly two lengths
\begin{eqnarray*}
L_1  &=  & | \langle \langle q_{n-1} \theta \rangle \rangle|  \\
L_2  &=&   |  \langle \langle q_{n-2}\theta \rangle \rangle  + 
(a_{n}-1) \langle \langle q_{n-1} \theta \rangle \rangle|,
\end{eqnarray*}
with $L_2 < L_1$. We deduce as in the case $a_{n+1} \ge 2$ that 
$$
L_2 < L_1 = L_2 + |\langle \langle q_n \theta \rangle \rangle| < 2L_2,
$$
as required.
$~~~\bsq$\\

The point of Lemma~\ref{lem21} is that for the choice $N=q_{n}-1$ the points
$\{ x+ j \theta ~(\bmod ~1): 0 \le j \le N\}$ are very close to uniformly spaced
on the interval $[0,1]$.  The next result obtains information on
the convergent denominators $q_{n}$ for the irrational number $\alpha_0$.

%
%

\begin{lemma}~\label{lem22}
For the irrational number  $\alpha_0= \log_3 2$ the
following hold. 

(1) For all  $q \ge 1$, and all integer $p$ there holds the
Diophantine inequality 
\beql{202}
| \alpha_0 - \frac{p}{q}| \ge \frac{1}{1200} \frac{1}{q^{c_0 + 1}}. 
\eeq
with $c_0=13.3$.

(2) The denominators $q_n$ of the continued fraction 
convergents $\frac{p_n}{q_n}$ of
$\alpha_0$ satisfy
\beql{220b}
q_n \le   1200  (q_{n-1})^{c_0}.
\eeq
\end{lemma}

\paragraph{Proof.} 
(i) The existence of
 a bound  of this general form, aside from the precise
 constants, follows from A. Baker's results on
linear forms in logarithms \cite[Theorem 3.1]{Ba75},
applied to  the linear form
$\Lambda = k + q\log 2 - p \log 3$, taking $k=0$, 
noting that its height $B := \max \{|p|, q \} \le 2q$.

The particular  bound \eqn{202} is obtained from a result of  
Simons and de Weger \cite[Lemma 12]{SW04}, who show that
for $k \ge 1$ and all integers $l$, 
$$
| (k +l)\log 2 - k \log 3| > \exp( -13.3(0.46057)) k^{-13.3}> \frac{1}{484} 
 k^{-13.3}.
$$
Their result is
proved using  a transcendence result of 
G. Rhin \cite[Proposition, p. 160]{Rh87}
for linear forms in two logarithms. We may suppose $k< k+l < 1.6k$, 
and obtain
$$
| \log_3 2 - \frac{k}{k+l}| > \frac{1}{\log 3}\exp( -13.3(0.46057)) 
(k+l)^{-1} k^{-13.3}
\ge \frac{1}{1200} (k+l)^{-14.3},
$$
which on taking $p=k, q=k+l$ gives the needed bound.\\

(2) Since $\alpha_0$ lies  in the interval between two successive
continued fraction convergents $\frac{p_{n-1}}{q_{n-1}}$ and $\frac{p_n}{q_n}$, 
we obtain using \eqn{202}  that 
$$
\frac{1}{q_n q_{n-1}} =| \frac{p_n}{q_n}- \frac{p_{n-1}}{q_{n-1}}|
= |\alpha_0 - \frac{p_{n-1}}{q_{n-1}}| + |{\alpha_0} - \frac{p_n}{q_n}|
\ge \frac{1}{1200}  \frac{1}{(q_{n-1})^{c_0+1}}
$$
Multiplying by $1200 q_n q_{n-1}^{c_0}$ gives \eqn{220b}.
$~~~\bsq$
%
%
\paragraph{Proof of Theorem~\ref{th11}.}
Let $\lambda >0$. We study for $1 \le n \le X$ the ternary expansion of
$$
x_n= x_n(\lambda)= \lfloor \lambda 2^n\rfloor.
$$
We will study the first $k$ leading ternary digits of  the $\{ x_n: 1 \le n \le X\}$
where we choose $k$ as follows. If $\frac{p_j}{q_j}$ are the convergents of
the continued fraction expansion of $\alpha_0= \log_3 2$, pick that $l$ such
that $q_{l-1} < X \le q_l$, and then choose $k$ to be the number of ternary
digits in $q_{l-1}$, so that $3^{k-1} < q_{l-1} \le 3^{k}$. Note that
$k = \lceil \log_3 q_{l-1} \rceil \le \lceil \log_3 X\rceil.$

We now set $w_n:= \log_3(\lambda 2^n) (mod~1),$ with $0 \le w_n <1,$ so
that
\beql{230}
w_n= n\alpha_0 + \log_3 \lambda ~(\bmod~1).
\eeq
We now observe that where $w_n$ falls in the interval $[0,1)$ specifies
the first $k$ ternary digits in the ternary expansion of $e^{w_n}$,  with
$1 \le e^{w_n} < 3$, we can  partition the interval $[0,1)$ into half-open intervals corresponding to
each such ternary expansion. Consider a  ternary expansion 
$$
\bb= [ b_0 b_1 \cdots b_{k-1}]_3,~~~ b_i \in \{0, 1,2\}, ~ b_0 \ne 0,
$$
of length $k$, noting there are $2 \cdot 3^{k-1}$ such expansions. Set
\beql{231}
\beta(\bb)= \sum_{j=0}^{k-1} \frac{b_j}{3^j},
\eeq
which has $1 \le \beta(\bb) < 3$ and associate the subinterval of $[0,1)$,
\beql{232}
J({\bb}) := [ \log_3 \beta(\bb), \log_3(\beta(\bb) + \frac{1}{3^{k-1}} ) ).
\eeq
These $2 \cdot 3^{k-1}$ subintervals partition $[0,1)$, from
$J([10\cdots 0]_3) =[ \log_3(1), \log_3 (1+ \frac{1}{3^{k-1} } ) )$ to
$J([22\cdots 2]_3)= [\log_3(3- \frac{1}{3^{k-1}}), \log_3 3).$

We claim that the following conditions (C1) and (C2)
are equivalent for $x_n$ with $3^m \le x_n \le 3^{m+1}$,
with $m \ge k$. 

(C1) $x_n$ has ternary expansion having the $k$ leading digits $\bb= [b_0 b_1 \cdots b_{k-1}]_3$, i.e
$x_n = \sum_{j=0}^m b_j 3^{m-j}$, for some $(b_{k+1},..., b_{m})$. 

(C2) $w_n= \log_3 x_n~ (\bmod ~1)$ has $w_n \in J(\bb)$.

\noindent The claim follows because the definition of $J(\bb)$ specifies the k leading ternary digits
of $3^{w_n}$, while $x_n= 3^m 3^{w_n}$ and the effect of multiplying by
$3^m$ simply shifts all ternary digits $m$ places to the left without changing the
leading digits. 

Next we note that the intervals $J(\bb)$ all have the same length to within
a factor of $3$, namely
\beql{235b}
\frac{1}{3^k} \le | J(\bb)|  \le \frac{1}{3^{k-1}}.
\eeq
This holds using 
$$|J(\bb)|=\log ( \beta(\bb) + \frac{1}{3^{k-1}}) - \log (\beta(\bb))=
 \int_{\beta(\bb)}^{\beta(\bb)+ \frac{1}{3^{k-1}}} \frac{dx}{x},
$$
and the bounds \eqn{235b} follow since 
 $\frac{1}{3} \le \frac{1}{x} \le 1$.

Next we examine the $w_n$ in consecutive blocks of length $N= q_{l-1}-1$,
i.e the set $\{ w_n:  j(q_{l-1}-1) \le n < (j+1)(q_{l-1}-1)\}.$ 
By \eqn{230} we may  apply  Lemma~\ref{lem21}(3) to this sequence of numbers,
to infer that the spacings between them are of two lengths
$L_1$ and $L_2$ which satisfy $L_2< L_1 < 2L_2$.  In particular since
$3^{k-1}  \le q_{l-1} \le 3^{k}$ these block sizes satisfy
$$
\frac{1}{2 \cdot 3^k} \le \frac{1}{2(q_{l-1} -1)} \le L_1 < L_2 \le \frac{2}{q_{l-1}-1} \le \frac{2}{3^{k-1}}.
$$
We conclude  using \eqn{235b} that at each subinterval $J(\bb)$ contains at most six
points $w_n$ from this block. Thus at most six values of $n$ in $j(q_{l-1}-1) \le n < (j+1)(q_{l-1}-1)$
give an $x_n$ having given intial $k$-digit ternary expansion $\bb=[b_0 b_1 \cdots b_{k_1}]_3$.

  We know there  are exactly $2^{k-1}$ values of $\bb=[b_0 b_1 \cdots b_{k_1}]_3$ that omit
  the ternary digit $2$, so the above shows there are at most $6 \cdot 2^{k-1}$
  values of $n$ in each such block giving an $x_n$ whose initial $k$ ternary digits
  avoid $2$. There are $\lfloor \frac{X}{q_{l-1}-1}\rfloor +1$ such blocks covering
  all $1 \le n \le X$ hence  we conclude there are at most
  \begin{eqnarray*}
 M := 
  6 \cdot 2^{k-1}\left(  \frac{X}{q_{l-1}-1} +1 \right) & \le & 6 \cdot 2^{k-1} \left( \frac{X}{3^{k-1}} + 1 \right)\\
  &\le & 6 \left( (\frac{2}{3})^{k-1} X + 2^{k-1}\right) \le 12 (\frac{2}{3})^{k-1} X,
  \end{eqnarray*}
values of $x_n$ whose initial $k$ ternary digits omit the digit $2$.
(In the last inequality we used $X \ge q_{l-1} > 3^{k-1}.$ 

It remains to upper bound $M$ as a function of $X$. Using Lemma~\ref{lem22}(2) 
we have
$$
X \le q_{l} \le 1200 (q_{l-1})^{c_0} \le 1200 (3^k)^{c_0}
$$
with $c_0=13.3$. We apply this bound to obtain
$$
(\frac{3}{2})^k = \left( 3^{c_0 k}\right) ^{\log_3(3/2) c_0^{-1}} \ge
 \left(\frac{1}{1200} X\right)^{(\frac{1-\alpha_0}{c_0})},
$$
Here $\frac{1}{37}< (\log_3(3/2)) c_0^{-1}= \frac{1- \alpha_0}{c_0} \le \frac{1}{36}$,
so we obtain
$$
(\frac{2}{3})^k \le (1200)^{\frac{1-\alpha_0}{c_0}} X^{- (\frac{1-\alpha_0}{c_0})}
$$
Substituting this into the definition of $M$ we obtain, 
$$
M \le 18 (\frac{2}{3})^k X \le 
18 \cdot (1200)^{\frac{1}{36}} X^{ 1- \frac{1-\alpha_0}{c_0}}
\le 25 X^{\frac{36}{37}} \le 25 X^{0.9725}.
$$
and the result follows. $~~~\bsq$.\\

%
%

\paragraph{Proof of Theorem~\ref{th12}.}
We will construct a rapidly
 increasing sequence of integers $S_0= \{m_k: k \ge 1\}$
having the form
\beql{281}
m_k = l_0+l_1+ ...+ l_k, 
\eeq
such that there is  an  uncountable
 set of real numbers $\tilde{\bbS}$ such that {\em all} the numbers $\lambda \in \bbS$
have the property:  for each $k \ge 1$,
 the integer $M_k := \lfloor \lambda 2^{m_{k}} \rfloor$  has
a  ternary expansion that omits the digit $1$. 
We now claim that all the integers  $N_k := \lfloor \lambda 2^{m_{k}-1} \rfloor$
have ternary expansions $(N_k)_3$ that omit the digit $2$. This holds 
because for each $N_k$  either $M_k = 2N_k$ or $M_k=2N_k+1$, but 
 $M_k$ is necessarily an even integer since all its
ternary digits are $0$ or $2$, so  we must have $M_k=2N_k$.
Thus $N_k$ has only digits $0$ and $1$ in its ternary expansion, so  we have
for $S=\{ m_k -1: k \ge 1\}$ that
$$
\tilde{\Sigma} \subset \Sigma (S):= \{ \lambda: ~(\lfloor \lambda 2^{n_k} \rfloor)_3 ~\mbox{omits~the~digit}~~2\},
$$
hence $\Sigma(S)$ is an uncountable set. 


 
We choose the $l_k$ recursively, taking $l_0=m_0=0$
and $l_k$  to be the smallest integer satisfying  $l_k \ge 2k$ and 
\beql{282}
0 < \{\{ \log_3 2^{l_k} \}\} = \{\{ l_k \alpha_0\}\} < 2^{-m_{k-1} - 2k -4}.
\eeq
Here $m_k= l_0+l_1+ \cdots +l_k$.
We set
$$
r_k := \lfloor l_k \alpha_0 \rfloor,   ~~~~~~\alpha_0= \log_3 2.
$$ 
The condition $l_k \ge 2k$ ensures that $r_k \ge k$.
Then we have 
$$
2^{l_k}= 3^{l_k \alpha_0}=3^{r_k + \{\{ l_k \alpha_0\}\}} =3^{r_k} 3^{\{\{ l_k \alpha_0\}\}}.
$$
Using $e^x \le 1+2x$ for $0 \le x \le 1$ we have
$$
3^{\{\{ l_k \alpha_0\}\}} = e^{\{\{ l_k \alpha_0\}\}\log 3} \le 1+ 2\log 3\{\{ l_k \alpha_0\}\}
\le 1+ \frac{2 \log 3}{2^{m_{k-1} + 2k + 4}}.
$$
Thus we obtain
\beql{283}
3^{r_k} < 2^{l_k} < 
3^{r_k}\left(1 + \frac{2 \ln 3}{2^{m_{k-1} + 2k +4}}\right)\le 
3^{r_k}\left( 1+ \frac{1}{3^{(m_{k-1} + 2k+2) \alpha_0}}\right)
\eeq
This says that the ternary expansion of $2^{l_k}$ has leading digit $1$ followed by
a string of at least $(m_{k-1} + 2k +2) \alpha_0$ zeros.

Given this choice of $\{ l_k: k \ge 1\}$, we define the set $\bbS$ to consist of all real numbers 
\beql{284}
\tilde{\Sigma} := \{ \lambda := \sum_{k=0}^{\infty} \frac{d_k}{2^{m_k}}: \lambda ~\mbox{is~ admissible}\}
\eeq
where $\lambda$ is called {\em admissible} if, for all $k \ge 1$ it has
the two properties

(P1) The digit $d_k$ satisfies
\beql{285}
0 \le d_k \le 3^{r_k} - 3^{r_k-k}.
\eeq

(P2) Let  $\lambda_k := \sum_{j=0}^k \frac{d_j}{2^{m_j}}$. Then the integer
\beql{286}
M_k := \lambda_k 2^{m_k}
\eeq
has a  ternary expansion $(M_k)_3$ which  omits the digit $1$. \\

{\bf Claim 1.} {\em Any $\lambda= \sum_{j=0}^{\infty} \frac{d_j}{2^{m_j}}$
with all $d_k$ satisfying (P1) satisfies
\beql{287}
1 \le \lambda < 2
\eeq
and
\beql{288}
M_k = \lambda_k 2^{m_k} = \lfloor \lambda 2^{m_k} \rfloor, ~\mbox{for~all}~~k \ge 1.
\eeq
}\\

To prove the claim , we observe that (P1) gives 
\begin{eqnarray}~\label{295}
1 \le \lambda &\le &
1 + \sum_{k=1}^{\infty} \frac{1}{2^{m_{k-1}}}
 \left( \frac{3^{r_k}- 3^{r_k -k}}{2^{l_k}}\right)  \nonumber\\
&\le & 1 + \sum_{k=1}^{\infty} 
 \frac{1}{2^{m_{k-1}}}(1- 3^{-k}) < 2. 
\end{eqnarray}
Next, (P1) gives 
\begin{eqnarray*}
0 \le \lambda - \lambda_k &=& \sum_{j=k+1}^{\infty} \frac{d_j}{2^{m_j}} 
=  
\frac{1}{2^{m_k}} \left( 
\sum_{j=k+1}^{\infty} \frac{d_j}{2^{m_{j} - m_{k}}}\right) \\
&\le & 
\frac{1}{2^{m_k}}\left( \sum_{j=k+1}^{\infty} (1 - \frac{1}{3^j})
\frac{1}{2^{m_{j-1}- m_{k}}}\right) \\
&\le & 
\frac{1}{2^{m_k}}\left( \sum_{j=k+1}^{\infty} (1 - \frac{1}{3^j})
\frac{1}{2^{(j-k-1)(2j)}} \right) <\frac{1}{2^{m_k}},
\end{eqnarray*}
proving Claim 1.\\

{\bf Claim 2.} {\em For any choice of $\{d_j: 1 \le j \le k-1\}$ that satisfy
both (P1) and (P2), there are at least $2^{r_k}- 2^{r_k-k}$ choices
of $d_k$ that satisfy (P1) and (P2).}\\

To prove this, first note that 
\beql{297}
\lambda_{k-1}2^{m_k} = M_{k-1} 2^{m_k - m_{k-1}}
= M_{k-1} 2^{l_k} = M_{k-1}3^{r_k} + M_{k-1}(2^{l_k}- 3^{r_k}).
\eeq
We assert that 
\beql{298}
0 \le  M_{k-1}(2^{l_k}- 3^{r_k}) \le 3^{r_k -k}.
\eeq
The left inequality is immediate,
and using \eqn{295} we have 
$M_{k-1} \le \lambda 2^{m_{k-1}} \le 2^{m_{k-1}+1}$,
while \eqn{283} gives
\begin{eqnarray*}
M_{k-1} (2^{l_k}- 3^{r_k}) &\le &
 2^{m_{k-1}+1} \left( 3^{r_k} \frac{\ln 3}{2^{m_{k-1} + 2k +4} } \right) \\
&\le&  3^{r_k} \frac{1}{2^{2k+3}} \le  3^{r_k-k},
\end{eqnarray*}
proving \eqn{298}.

From \eqn{297} and \eqn{298} we see that the ternary expansion
of $\lambda_{k-1}2^{m_k}$ repeats that of $M_{k-1}$ shifted
$r_k$ positions to the left, then has a block of at least
$k$ zeros, and following this  has the ternary expansion of  the integer
$ M_{k-1}(2^{l_k}- 3^{r_k})$. It follows that choosing from the range of values 
 $0 \le d_k \le 3^{r_k} - 3^{r_k-k}$, and setting $\lambda_k := \sum_{j-0}^k \frac{d_j}{2^{m_j}}$,
 the  integers 
\beql{299}
M_{k} := \lambda_{k}2^{m_k}= \lambda_{k-1}2^{m_k} + d_k
\eeq
can be selected to give all ternary integers which

(i) have the ternary expansion matching $M_{k-1}$ to the left of the
$r_{k}$-th position, 

(ii) omit the digit $1$, and

 (iii) have at least one $2$ and at least one $0$ in
positions between $r_k$ and $r_k -k$;

\noindent  call these {\em allowable values.}  In these $k+1$
positions the largest allowed value is $222\cdots 20$
and the smallest is $000\cdots02$.  These produce exactly 
$2^{r_k} - 2^{r_k-k}$ such  ternary integers $M_k$, constructed by choice of
the same number of allowable values $d_k$. This proves Claim 2.\\

{\bf Claim 3.} {\em The set $\tilde{\Sigma}$ contains  uncountably many admissible $\lambda$,
and each of them has the property that every
\beql{299a}
M_k = \lfloor \lambda 2^{m_k} \rfloor, ~~~k \ge 1,
\eeq
has a ternary expansion $(M_k)_3$ that omits the digit $1$.}\\

Indeed Claim 2 implies there are uncountably many such $\lambda$,  since
the construction has a Cantor set form which
gives an infinite tree of values with branching at least two
at every node at  every level $k \ge 2$. The relation \eqn{299a} holds
by Claim 1, and these $M_k$  have ternary expansions omitting 2 by
(P2). Thus Claim 3
follows. 

It remains to verify the upper and lower bounds \eqn{103} on the 
growth rate of the sequence $m_k$.
The size of $m_k$ is determined by the Diophantine condition on $l_k$
given by equation \eqn{282}. 
(The numbers $l_k$ grow so rapidly that the side condition $l_k \ge 2k$
is automatically satisfied for $k \ge 2$.)   
Note that we cannot directly use Dirichlet's box principle to get an
upper bound for the size of the minimal $l_k$ satisfying  \eqn{282}
because this  is a one-sided approximation condition.
Instead we have that the minimal $l_k$ 
 will be no larger than that  even-numbered convergent $q_{2l}$
of the continued fraction expansion of $\alpha_0$ satisfying
$$ 
q_{2l-2} \le 2^{m_{k-1} + 2k+4} < q_{2l}.
$$
Lemma~\ref{lem22} (2) gives the bound 
\beql{299b}
q_{2l} \le \frac{1}{C_0^2} (q_{2l-2})^{2c_1} = (1200)^2(q_{2l-2})^{26.6}
\le 2^{27m_{k-1} + 54k + 132}.
\eeq
Since $n_k= m_k -1$ we obtain 
$$
n_{k} \le m_{k} \le m_{k-1} + q_{2l} \le m_{k-1} +  2^{27m_{k-1} + 54k + 132} \le 
 2^{27(n_{k-1}+ 2k+6)},
$$
which is the upper bound in \eqn{103}. 

Lemma~\ref{lem22}  implies a  lower bound on how
small $l_{k+1}$ can be to make \eqn{282} hold,
namely we must have 
\beql{249c}
(l_{k+1})^{c_0} \ge 2^{m_k +2j - 7},
\eeq
with $c_0=13.3,$ to avoid contradicting \ref{lem22}(1).
This yields the lower bound in \eqn{103}, which holds
for $n_k=m_k-1$ produced in this construction. 
$~~~\bsq$

%
%

\paragraph{Proof of Theorem~\ref{th13}.}
We consider the truncated exceptional 
set $\sE_{T}(\RR_{+}) $ . We 
first  establish the upper bound $dim_{H}( \sE_{T}(\RR_{+})) \le \alpha_0$.
We have
$$
\sE_{T}(\RR_{+}) = \bigcup_{M=2}^{\infty} \left(  \sE_{T}(\RR_{+})\cap [ \frac{1}{M}, M]\right).
$$
Since the Hausdorff dimension of a countable union of sets is the supremum of the
Hausdorff dimensions of the separate sets, it suffices to show that
\beql{260}
dim_{H}(  \sE_{T}(\RR_{+})\cap [ \frac{1}{M}, M]) \le \alpha_0 = \log_3 2.
\eeq
To show this we find suitable coverings of these sets. For each $n \ge 1$ we have
\beql{261}
\sE_{T}(\RR_{+})\cap [ \frac{1}{M}, M]) \subset S_n(M) := \bigcup_{j=N}^{\infty} \Sigma_j([\frac{1}{M}, M])
\eeq
with
$$
\Sigma_{j}([\frac{1}{M}, M]) := \{ \lambda: -\frac{1}{M} \le \lambda \le M~~\mbox{and}~
(\lfloor \lambda 2^j\rfloor )_3 ~\mbox{omits~the~digit}~2\}.
$$
The set $S_n(M)$ thus encodes a "tail event" that there are arbitrarily large 
$j$ for which $(\lfloor \lambda 2^j\rfloor )_3$ that omit the digit $2$. 
We will eventually let $n \to \infty$ so we suppose that $n \ge \log_3 M +2$, so that
$\lambda 2^j \ge 1$, for any $j \ge n$. Now consider such $j$ as fixed,
and note that $\lfloor \lambda 2^j \rfloor$ takes a fixed integer value on
an interval of length $\frac{1}{2^j}$. Letting $\bb = (\lfloor \lambda 2^j\rfloor )_3$,
we see that allowable values of $\bb$ satisfy $1 \le \bb \le M 2^j$. As $\lambda$
varies over $[\frac{1}{M}, M]$ these integers vary over a subset of $[1, M 2^j]$ and
of these, the number of such ternary expansions $\bb$ that omit the digit $2$ is
at most (counting integers over successive blocks $[3^{k-1}, 3^k)$), 
\begin{eqnarray*}
1+ 2 + \cdots+ 2^{ \lceil \log_3 (2^j M)\rceil} & \le & 2^{ \log_2(2^j M) + 2}\\
&\le & 2^{j \alpha_0 + \log_3 M + 2} \le 4M 2^{j \alpha_0}.
\end{eqnarray*}
Thus we obtain a collection 
$$\sI_j(M):= \{ I_j(\bb): ~\bb ~~\mbox{gives~an~admissible~interval~for}~ \lfloor \lambda 2^j\rfloor, ~\frac{1}{M} \le \lambda \le M\}.$$
of at most $4M 2^{j \alpha_0}$ intervals of length $\frac{1}{3^j}$, and these intervals
cover the set 
$\Sigma_{j} ([\frac{1}{M}, M] ).$
Summing over all
$j \ge n$ we obtain an infinite collection of intervals 
$$
\sI(n, M) :=\bigcup_{j =n}^{\infty} \sI_j(M),
$$
which  cover the set $\sE_{T}(\RR_{+})\cap [ \frac{1}{M}, M])$ by \eqn{261}, 
and every interval included has length at most $\frac{1}{2^n}.$ 
Now fix $\epsilon >0$
and observe that
\begin{eqnarray*}
\sum_{ I \in \sI(n ,  M)} |I|^{\alpha_0 + \epsilon} & =& 
\sum_{j=n}^{\infty} \left( \sum_{I \in \sI_j(M)} (\frac{1}{2^j})^{\alpha_0 + \epsilon} \right) \\
&\le & \sum_{j=n}^{\infty} 4M 2^{j \alpha_0} (\frac{1}{2^j})^{\alpha_0 + \epsilon} \\
&=& 4M \left( \sum_{j=n}^{\infty} 2^{- j \epsilon}\right) 
= (\frac{4M}{1- 2^{- \epsilon}} ) 2^{- n \epsilon}.\\
\end{eqnarray*}
Letting $n \to \infty$, the diameter of the covering  $\sI(n, M)$   goes to zero, and the scaled length goes
to zero as well, which establishes
$$
dim_{H} \left( \sE_{T}(\RR_{+})\cap [ \frac{1}{M}, M] \right)\le \alpha_0 + \epsilon.
$$
Now we can let $\epsilon \to 0$ to obtain \eqn{260}, and the upper bound 
$dim_{H}( \sE_{T}(\RR_{+})) \le \alpha_0$ follows. 
%
%

To  establish  the lower bound $dim_{H}( \sE_{T}(\RR)) \ge \alpha_0$
is more difficult, as it requires controlling all coverings of the set.
We  will actually establish  the stronger result that
\beql{262}
meas_{\alpha_0}( \tilde{\bbS}) >\frac{1}{16},
\eeq
where $\tilde{\Sigma} \subset [1,2]$ is the set constructed in
Theorem ~\ref{th12} in \eqn{284}.
The set $\tilde{\Sigma}$ had a  construction resembling a Cantor set,
with two differences. The first difference is that
 the dissection at each layer $k$ depended on the previous layers,
 and the second difference 
 is that the layer at level $k$ involved denominators $2^{m_k}$ with
$$
m_k = l_0+l_1+ ...+ l_k, 
$$
with the $l_k$ growing extremely rapidly.
We can however  adapt an argument given
in Falconer \cite[Example 2.7, p. 31]{Fa90} for the 
Cantor set to show \eqn{262}.

We claim  that $\tilde{\bbS}$ has a representation as 
\beql{263}
\tilde{\bbS} = \bigcap_{s=1}^{\infty} X_s,
\eeq
in which $X_s$ consists of a union of a collection $\sJ_s$ of disjoint intervals
of size proportional to $3^{-s}$, and the sets are nested: 
$$
\cdots X_3 \subset X_{2} \subset X_{1}.
$$ 
Here the intervals in $\sJ_s$ will play the role of the Cantor set dissection into intervals
at level $s$,  for each power of $3^s$.

We first define the collection $\sJ_s$ for those levels $s=s_k$ with
\beql{264} 
s_j: = \lfloor m_j \alpha_0 \rfloor,
\eeq
which are directly given
in the construction of Theorem~\ref{th12}. Then we
show one can fill in all the intermediate layers $s_k \le s < s_{k+1}$.  

We have
$3^{s_k} < 2^{m_k} < 3^{s_k +1}$, and the set $\sJ_{s_k}$ is the
union of all closed intervals
$$
\sJ_{s_k}:=\{ \left[\frac{M}{2^{m_k}},\frac{M+1}{2^{m_k}}\right]: ~M= \lambda_k 2^{m_k}
~\mbox{with}~\lambda_k = \sum_{j=0}^k \frac{d_j}{2^{m_j}}~\mbox{admissible}.\}
$$
with admissibility in the  construction in Theorem~\ref{th12}.
Here we have
$$
2^{m_k} = 2^{l_1 +...+l_k} = 3^{l_1\alpha_0+...+l_k\alpha_0} 
=3^{r_1+r_2+...+r_k} \cdot 
3^{ \{\{ l_1 \alpha_0\}\} + ... + \{\{ l_k \alpha_0\}\} }
\le  2 \cdot  3^{r+1 +... +r_k}, $$
using the fact that 
$$
\sum_{k=1}^{\infty} \{\{ l_k \alpha_0\}\} \le \sum_{k=1}^{\infty} 2^{-m_{k-1}-2k-2} \le \frac{1}{2},
$$
using \eqn{282}. This also establishes that
\beql{265}
s_k= r_1 +r_2+...+r_k.
\eeq
Inside each interval at level $s=s_{k-1}$ there fit exactly
$2^{r_k}- 2^{r_k-k}$ subintervals at ternary level  $s=s_k$,
each of length $2^{-m_k}$, and we now know that
$\frac{1}{2} 3^{-s_k} \le 2^{-m_k} \le 3^{-s_k}.$
This dissection of an interval at ternary level $s_{k-1}$
into subintervals at  ternary level $s_k$ 
is exactly that of the Cantor set,
except that the two  ends of the interval are trimmed off
a small amount, to a relative distance  $  3^{-k}$ from
each end of the interval. 

We now  fill in the intermediate
levels $X_{s}$ for $s_{k-1} < s <  s_{k}$ by gluing
together all intervals  in $\sJ_{s_k}$ that have matching initial ternary expansions
$[M]_3$ of $M=\lambda_k 2^{m_k}$, disregarding the
last $s_{k}-s$ ternary digits of $[M]_3$, and
filling in the space between them.  The resulting
intervals of $\sJ_{s}$ all have size exactly $3^{s_k-s} 2^{-m_k}$
(except possibly for two subintervals adjacent to the truncated ends);
their size lies between $\frac{1}{2} 3^{-s}$ and $3^{-s}$. Also, 
 the gaps between any two adjacent intervals at ternary level $s$
are of size at least as large as 
\beql{266}  
G_s = 3^{s_k-s} 2^{-m_k} \ge \frac{1}{2} 3^{-s}.
\eeq
This  fact holds because this construction uses ternary integers
omitting the digit $1$; the set of  ternary integers
omitting the digit $2$ has some intervals of this kind that 
 are adjacent, so the gap size would be zero in that case.

The above construction defines the intervals in  $\sJ_s$ at level $s$ for all $s$.
This dissection imitates the Cantor set in that each interval
at level $s$,  contains at most $2^{s'-s}$ subintervals at any deeper
ternary level $s' \ge s$. It may contain fewer subintervals, due to
the trimming at ends of the subinterval, but it always 
contains at least $2^{s'-s-1}$ such subintervals.

The set $\tilde{\bbS}$ is a compact set contained in
the interval $[1,2]$. To bound its $\alpha_0$-dimensional Hausdorff measure 
from below,
we must show that in {\em every}  covering $\{U_i\}$ by closed intervals
there holds
\beql{267}
\sum_{i} |U_i|^{\alpha_0} \ge \frac{1}{16}.
\eeq
By enlarging the intervals slightly (by $1+\epsilon$) and
observing that their interiors give an open cover of $\tilde{\bbS}$,
we can extract a finite subcover.  
Since we can extract a finite subcover for any $\epsilon>0$,  it suffices
to verify  \eqn{267} holds for every finite cover $\{U_i\}$ of $\tilde{\bbS}$ 
by intervals.

Given an interval $U_i$ in a covering, define $s$ by 
\beql{268}
3^{-s} \le |U_i| < 3^{-s+1}.
\eeq
Then $U_i$ can touch at most two subintervals at level $s$
because all subintervals in $\sJ_{s}$ are sepated by gaps of size at least  
$\frac{1}{2} 3^{-s}.$ If $s' \ge s$ then $U_i$ intersects
at most $2 \cdot 2^{s'-s}$ subintervals at level $s'-s$; 
by \eqn{268} this number  is bounded above by
\beql{268a}
2 \cdot  2^{s'-s} \le 2^{s'} 3^{-\alpha_0 s} \le 
2 \cdot 2^{s'} (3^{\alpha_0} |U_i|^{\alpha_0}) = 
4  \cdot  2^{s'}|U_i|^{\alpha_0}).
\eeq
Given a finite cover,  choose $s'=s_k$ large enough so that 
$|U_i| \ge 3^{-s'}$ for all $i$. Then the collection   $\{U_i\}$
necessarily covers all subintervals at level $s'= s_k$. By 
 construction $\sI_{s_k}$ contains at least 
\beql{269}
\prod_{i=1}^k ( 2^{r_i} - 2^{r_i -i}) = 
2^{r_1+...+r_k} \prod_{i=1}^n(1-2^{-i}) \ge \frac{1}{4}  2^{s_k} 
\eeq
intervals, since 
where  $\prod_{i=1}^k(1-2^{-i}) \ge  \prod_{i=1}^{\infty}(1-2^{-i}) \ge \frac{1}{4}$.
Now we count how many intervals at level $s_k$ are covered. Since
$U_i$ intersects at most $4 \cdot 2^{s_k}|U_i|^{\alpha_0}$ such intervals we
must have
$$
\sum_i 4 \cdot 2^{s_k}|U_i|^{\alpha_0} \ge |\sJ_{s_k}| \ge \frac{1}{4} 2^{-s_k}.
$$
This yields 
$$
\sum_{i}|U_i|^{\alpha_0} \ge \frac{1}{16},
$$
which establishes \eqn{262}.
$~~~\bsq$

\paragraph{Remark.}
More  generally we may consider the  real dynamical system $y \to \beta y$, where $\beta>1$,
and consider the truncated ternary expansions  $\{(\lfloor \lambda \beta^n \rfloor)_3: n \ge 0\}$. 
The methods above should extend to those $\beta$ such that
$\alpha:= \log_3 \beta$  satisfies a Diophantine condition
\beql{270}
| \alpha - \frac{p}{q}| \ge c_2 \frac{1}{q^{c_1 + 1}}, ~~\mbox{for~all} ~p,q ~\mbox{with} ~q \ge 1, 
\eeq
for  constants $c_1>1$ and $c_2>0$. The conclusions of the results require appropriate
modification, with constants depending on the Diophantine condition. 

%
%
%
%
\section{$3$-adic Integer Dynamical System: Proofs}

We consider the $3$-adic integers $\ZZ_3$ and write the $3$-adic
expansion of $\lambda \in \ZZ_3$ as
\beql{310}
\lambda= \sum_{j=0}^{\infty} d_j 3^j~~~\mbox{with~each}~d_j \in \{ 0, 1, 2\}.
\eeq
We write the $3$-adic digit expansion as ~$(\lambda)_3 = ( \cdots d_2 d_1 d_0)_3.$

This dynamical system  consider the sequence of $3$-adic integers,
$y_n = \lambda 2^n,$ where $\lambda$ is a given
nonzero $3$-adic integer. Here $y_n$
form the forward orbit of  the first order linear
recurrence
$y_n = 2y_{n-1}$, with initial condition
$y_0= \lambda$.   The map $T: x \to 2x$ is
an automorphism of  the $3$-adic integers $\ZZ_3$, which leaves  each of the 
sets  $\Sigma_j :=  3^j \ZZ_3^{\ast}$ for $j \ge 0$  invariant.
(Here $\ZZ_3^{\ast}$ are the $3$-adic units.) These sets partition $\ZZ_3$ and 
this map acts  ergodically on each  component $\Sigma_j$.

We are interested in the possible ways
that the orbit $\{ y_n: n \ge 0\}$ can
intersect the set
$
\Sigma_{3, \bar{2}} := \{w: w = \sum_{j=0}^{\infty} a_j 3^j \in \ZZ_3, 
~\mbox{with~each}~ a_j= 0 ~\mbox{or}~ 1 \}.
$
We now  upper bound the number of $n \le X$ that
can fall in the set $\Sigma_{3, \bar{2}}$.

%
%
\paragraph{Proof of Theorem~\ref{th14}.}
Let $\lambda \in \ZZ_3$ with $\lambda \ne 0$. We study the set
\beql{311}
\tilde{N}_{\lambda}(X) := \# \{1\le n \le X:~ (\lambda 2^n)_3~~\mbox{omits~the~digit}~2\}.
\eeq
Write $\lambda= 3^j \lambda^{\ast}$ with
$\lambda^{\ast} \in \ZZ_3^{\times}:= \{ \lambda \in \ZZ_3: ~\lambda \not\equiv 0~(\bmod~3)\}.$
Then we have $\tilde{N}_{\lambda}(X) = \tilde{N}_{\lambda^{\ast}}(X)$, since multiplication
by $3^j$ simply shifts $3$-adic digits to the left. Thus to prove the desired inequality there
is no loss of generality to require
$\lambda \ne 0~(\bmod~3)$, by replacing $\lambda $ with $\lambda^{\ast}$. 

The proof is based on the fact that $2$ is
a primitive root $(\bmod~3^k)$ for each $k \ge 1$. Thus, for each $k \ge 1$
\beql{313}
\{ \lambda 2^n~(\bmod~3): 1 \le n \le \phi(3^k)= 2 \cdot 3^{k-1} \}
\eeq
runs over all $2 \cdot 3^{k-1}$ invertible residue classes $(\bmod~3^k)$. Of these,
exactly $2^{k-1}$ residue classes have a $3$-adic expansion that omits the digit $2$.
Now, given $X$, pick that $k$ such that
$$
2 \cdot 3^{k-2} < X \le 2 \cdot 3^{k-1}.
$$
Then  applying \eqn{313} over $1 \le n \le 2 \cdot 3^{k-1}$ we have exactly $2^{k-1}$
values of $n$ with $(\lambda 2^n)_3$ omitting the digit $2$ in its first $k$
$3$-adic digits $(d_{k-1} \cdots d_1 d_0)_3.$  Thus
\begin{eqnarray*}
\tilde{N}_{\lambda}(X) & \le & 2^{k-1} = 2 \cdot 2^{k-2} = 2 \cdot 3^{\alpha_0(k-2)}\\
&= & 2^{1-\alpha_0} \left (2 \cdot 3^{k-2}\right)^{\alpha_0} \le  2 X^{\alpha_0}, 
\end{eqnarray*}
which is the desired upper bound. $~~~\bsq$ \\

The object of Theorem~\ref{th15} is to establish
upper bounds on the  Hausdorff dimension of the 3-adic exceptional
set $\sE(\ZZ_3)$ through upper bounds on various
$\sE^{(j)}(\ZZ_3)$ which contain it. 
  
 We note that 
 Hausdorff dimension is a metric notion (cf. Rogers \cite{Ro70}),
and its version for $3$-adic integers uses the $3$-adic metric is 
quite similar 
to Hausdorff dimension for real numbers on the interval $[0,1]$. In fact we
have a continuous (and almost one-to-one) mapping
$\iota: \ZZ_3 \to [0,1]$ which sends a $3$-adic number $\lambda= (  \cdots d_2d_1d_0)_3$ 
to the real number with ternary expansion $.d_0d_1d_2 \cdots$.
One can show that this mapping preserves Hausdorff dimension of sets, 
i.e a $3$-adic set $X$ and its image $\iota(X)$ have the same Hausdorff dimension.
This holds because one can expand each set in  a  $3$-adic covering  of a set $X$
to a  closed-open disk\\
$B(m, 3^j) = \{ x \in \ZZ_3: ~x \equiv m~(\bmod~3^j)\}$, with at most a factor of $3$ increase in
diameter, and similarly one can inflate any real covering to a covering
with ternary intervals $[\frac{m}{3^j}, \frac{m+1}{3^j}]$ with at most a factor of
3 increase in diameter. But these special intervals are assigned the same diameter under
their respective metrics, and this can be used to show the Hausdorff dimensions
of $X$ and $\iota(X)$ coincide. 
 In particular the standard $3$-adic Cantor set $\Sigma_{3, \bar{1}}$
maps under $\iota$ to the usual Cantor set in $[0,1]$ hence it has Hausdorff dimension
$d_H(\Sigma_{3, \bar{1}})= \log_3( 2) \approx 0.63092$.
Now $\Sigma_{3, \bar{1}}= 2 \Sigma_{3, \bar{2}})$ hence 
$\dim_{H}(\Sigma_{3, \bar{2}}) = \log_3(2)$ as well.

%
%
\paragraph{Proof of Theorem~\ref{th15}.}  This proof  assumes that
Theorem~\ref{th16} is proved in order to deduce the upper bound in (2).\\

(1) We have 
$$
\sE^{(1)}(\ZZ_3) = \bigcup_{m=0}^{\infty} \sC(2^m),
$$
with
$
\sC(2^m) :=\{ \lambda:~(\lambda 2^n)_3 ~\mbox{omits~the~digit}~~2\}.
$
Then
$$
\sC(2^m) = \frac{1}{2^m} \sC(1) = \frac{1}{2^m}( \Sigma_{3, \bar{2}}) =
 \frac{1}{2^{m+1}}( \Sigma_{3, \bar{1}}).
$$
Each $\sC(2^m)$ is a linearly rescaled version of the Cantor set 
$\Sigma_{3, \bar{1}}$
so has Hausdorff dimension $\log_3 2$. Thus
$$
\log_3 2 = \dim_{H}(\sC(1))\le \dim_{H}(\sE^{(1)}(\ZZ_3)) \le
 \sup_{m \ge 0} \dim_{H}(\sC(2^m)) = \log_3 2,
 $$
 as required.
\smallskip

(2) We have
$$
\sE^{(2)}(\ZZ_3) = \bigcup_{0 \le m_1 < m_2} \sC(2^{m_1}, 2^{m_2}).
$$
with 
$\sC(2^{m_1}, 2^{m_2}) :=\{ \lambda:~(\lambda 2^{m_i})_3 ~\mbox{omits~the~digit}~~2\}.$
Now
$$
\sC(2^{m_1}, 2^{m_2}) = \frac{1}{2^{m_1}} \sC(1, 2^{m_2- m_1}),
$$
which gives $\dim_{H} (\sC(2^{m_1}, 2^{m_2})) = \dim_{H} (\sC(1, 2^{m_2-m_1})).$
Since $m_2 - m_1 \ge 1$, Theorem~\ref{th16} applies to give
$$
\dim_{H} (\sC(1, 2^{m_2-m_1}) )\le \frac{1}{2}, ~~\mbox{for~all}~ m_2> m_1 \ge 0.
$$
This yields the upper bound
$$
\dim_{H}(\sE^{(2)}(\ZZ_3))= \sup_{0 \le m_1< m_2} \dim_{H}(\sC(2^{m_1}, 2^{m_2})) \le \frac{1}{2}.
$$
To establish the lower bound, we use the fact that $4=(11)_3$. Then the set 
$$
\Sigma_A:= \{ \lambda=(\cdots d_2d_1d_0)_3: ~ \mbox{all~blocks}~ d_{2n+1}d_{2n} \in \{
 00,01 \} ~\} \subset \Sigma_{3, \bar{2}},
$$
satisfies
$$
4\Sigma_A = \{ \lambda=(\cdots d_2d_1d_0)_3: ~ \mbox{all~blocks}~ d_{2n+1}d_{2n} \in \{ 00, 11\}
~\} \subset \Sigma_{3, \bar{2}},
$$
which shows that $\Sigma_A \subset \sC(1,4).$
Now $\Sigma_A$ is given by a Cantor set construction, which permits its Hausdorff dimension
to be computed in a standard way. We obtain
$$
\dim_{H}(\sE^{(2)}(\ZZ_3)) \ge \dim_{H}(\sC(1, 2^2)) \ge \dim_{H}(\Sigma_A) = \frac{\log_3 (2)}{\log_3(9)}
= \frac{1}{2} \log_3 (2) \approx 0.31596.
$$
\smallskip

(3) We have
$$
\sE^{(2)}(\ZZ_3) = \bigcup_{0 \le  m_1 < m_2< m_3} \sC(2^{m_1}, 2^{m_2}, 2^{m_3}).
$$
 The upper bound $\dim_{H} (\sE^{(3)}(\ZZ_3) \le \dim_{H}(\sE^{(2)}(\ZZ_3)$ is immediate.
 To establish the lower bound, we use the facts  that $4=(11)_3$ and
 $256=(100111)_3$. Then
$$
\Sigma_{B} := \{ \lambda=(\cdots d_2d_1d_0)_3: ~ \mbox{all}~ d_{6n+5}d_{6n+4}d_{6n+3}d_{6n+2}d_{6n+1}d_{6n} \in \{000000, 000001 \}~ \} \subset \Sigma_{3, \bar{2}}.
$$
has 
$$
4\Sigma_{B}=\{ \lambda=(\cdots d_2d_1d_0)_3: ~ \mbox{all}~
 d_{6n+5}d_{6n+4}d_{6n+3}d_{6n+2}d_{6n+1}d_{6n} \in \{ 000000,~000011 \}~ \} 
 \subset \Sigma_{3, \bar{2}}.
$$
$$
256\Sigma_{B}=\{ \lambda=(\cdots d_2d_1d_0)_3:  \mbox{all}~ 
d_{6n+5}d_{6n+4}d_{6n+3}d_{6n+2}d_{6n+1}d_{6n} \in \{ 000000, 100111\}~\}
 \subset \Sigma_{3, \bar{2}}.
$$
Thus $\Sigma_B \subset \sC(1, 4, 256) \subset \sE^{(3)}(\ZZ_3)$. Now $\Sigma_B$
has a Cantor set construction showing that
$$
\dim_{H}(\Sigma_B) = \frac{\log_3(2)}{\log_3 (3^6)} = \frac{1}{6} \log_3 (2) \approx 0. 10515,
$$
which gives the asserted lower bound. 
$~~~\bsq$

\paragraph{Remark.} The proof of Theorem~\ref{th15}
exploited the known solutions to Erd\H{o}s's problem.
Consequently this approach does not extend to give a
nonzero lower bound for $\dim_H(\sE^{(k)}(\ZZ_3))$,
for any $k \ge 4$.  Theorem~\ref{th17} offers more flexibility in finding 
ternary expansion identities for integers that could potentially yield nonzero lower bounds in
these cases.

%
%
%
%
\section{Intersections of Multiplicative Translates of the $3$-Adic Cantor Set: Proofs}

We study the $3$-adic Cantor set $\Sigma_{3, \bar{1}}$, defined by
\beql{501}
\Sigma_{3, \bar{2}} := \{ \lambda \in \ZZ_3:~\mbox{the~3-adic~digit~expansion}~(\lambda)_3~
\mbox{omits~the~digit~2} \}.
\eeq
For integers $1 \le M_1< M_2 < \cdots < M_k$ we  define the {\em intersection set}
\begin{eqnarray} \label{502}
\sC(M_1, M_2, \cdots , M_k) & := & 
\{ \lambda \in \ZZ_3:~ (M_i \lambda)_3 ~\mbox{omits~the~digit}~~2\} \\
&= &\bigcap_{i=1}^k  \frac{1}{M_i }\Sigma_{3, \bar{1}}
\end{eqnarray}
In  \S3 we used integers  $M_i = 2^{m_i}$ but here we
allow arbitrary positive integers $M_i$.  We study $\sC(1, M)$ for general $M$ and note first
that $\sC (1, 3^j M) = \sC(1, M).$. Thus without loss of generality
we may reduce to the case $gcd(M, 3)=1$. Another simple fact is the following. 

%
%

\begin{lemma}~\label{le41}
Let $M $ be a positive integer.\\ 

(1) If $M \equiv 2 (\bmod~3)$ then $\sC (1, M) = \{ 0\}$.\\

(2) If $M \equiv 1 (\bmod~3)$ then $\sC(1, M)$ is an infinite set. 
\end{lemma}

\paragraph{Proof.}
(1) Suppose $M \equiv 2 (\bmod~3)$.
If  $\sC(1, M) \ne \{0\}$, then it necessarily contains some $\lambda$ with
$\lambda \ne 0 (\bmod~3)$, since we may divide out any powers of $3$, and 
 multiplication by $3^j$ simply shifts digits to the left. Then $\lambda \in \Sigma_{3 , \bar{2}}$
implies $\lambda \equiv 1~(\bmod~3)$.  Then $M \lambda \equiv 2(\bmod~3)$ so 
$M\lambda \not\in \Sigma_{3, \bar{2}}$,  a contradicting membership
in $\sc(1,M)$. Hence no such $\lambda $ exist, and 
 $\sC(1, M) = \{0\}$.

(2) Suppose $M \equiv 1 (\bmod~3).$ To show $\sC(1, M)$ is an infinite set it 
suffices to exhibit one nonzero element $\lambda \in \sC^{\ast}(1, M)$, because
$3^j \lambda \in \sC^{\ast}(1, M)$ for all $j \ge 0$. We may construct such an element
$\lambda= ( \cdots d_2 d_1 d_0)_3$ recursively, starting with the choice $d_0=1$. 
Write $M = \sum_{j=0}^n a_j 3^j$, with $a_0=1$.
Let $M \lambda= \sum_{j=0}^{\infty} c_j 3^j$. Then the $k$-th digit satisfies
$$
c_k \equiv d_k + \left(\sum_{j=1}^{n} a_j d_{n-j} \right)+ e_{k-1} ~(\bmod~3)
$$
(with the convention $d_{-1} = d_{-2}= \cdots = d_{-n}=0$), and with
$e_{k-1}$ encoding the "carry digit" information, from the previous terms, which
is completely determined by $(d_0, d_1, ..., d_{k-1}.)$ Since we have two 
choices $0$, $1$ for $d_k$, at least one of them will foce $c_k \ne 2 ~(\bmod~3).$
Thus we can recursively construct an admissible $\lambda$ by
induction on $k$.  $~~~\bsq$.\\

It is possible to make a  detailed analysis of the structure
of $\sC(1, M)$ with $M \equiv 1 ~(\bmod~3)$, and determine their Hausdorff dimensions, 
which we consider elsewhere. One can show that 
 infinite set $\sC(1, M)$ can be either countable or uncountable, e.g. $\sC(1, 49)$ is
 countably infinite, while $\sC(1, 7)$ is uncountable.

Now we upper bound the Hausdorff dimension of $\sC(1, M)$. For $M=3^j,~(j \ge 0)$ we have
 $\sC(1, 3^j) = \Sigma_{3, \bar{2}} $, whence $\dim_{H}(\sC(1, 3^j)) = \log_3(2) \approx 0.63$.
 The following result treats all other $M \ge 1$.

%
%
\paragraph{Proof of Theorem~\ref{th16}.}
We suppose that $M>1$ is an integer that is not a power of $3$, i.e. its ternary
expansion $(M)_3$ contains at least two nonzero ternary digits.
Our object is  to upper bound the Hausdorff dimension of 
$$
\sC(1, M) := \Sigma_{3, \bar{2}} \cap M \Sigma_{3, \bar{2}},
$$
by $\frac{1}{2}$.
By the discussion above we may reduce to the case that $gcd(M, 3)=1$,
and by Lemma~\ref{le41} we may suppose $M \equiv 1~(\bmod~3),$ 
since the Hausdorff dimension is $0$ if $M \equiv 2 (\bmod~3).$ 
Thus we can write
\beql{510}
(M)_3= b_0 + b_m3^m + \sum_{j=m+1}^{n} b_j 3^j,~~~b_j \in \{0, 1,2\}, ~\mbox{with}~~ b_0b_m \ne 0. 
\eeq 
and  $b_0=1$, where the $m$-th digit is the first nonzero ternary digit after the $0$-th digit. 

We will study the 
minimal  covers of  $\sC(1, M)$ with $3$-adic open sets of measure $3^{-r-1}$ that specify
the first $r+1$ digits of the $3$-adic expansion of a number $\lambda \in \sC(1, M)$.
These sets are congruence classes $(\bmod~3^{r+1})$ and they have diameter
$3^{-(r+1)}.$ We call a congruence class $ \lambda ~(\bmod~3^{r+1})$ {\em admissible}
if $\sC^{\ast}(1,M)$ contains at least one element in this congruence class. 
Our object is to bound above the number of admissible congruence classes
$\lambda ~(\bmod ~3^{r+1})$
 
Set  $\lambda= \sum_{j=0}^\infty d_j 3^j \in \Sigma_{3, \bar{2}}$, so that 
 each $d_j =0$ or $1$. Now define the
digits $a_j$ by
$$
M \lambda = \sum_{j=0}^{\infty} a_j 3^j, ~~~a_j \in \{0, 1,2\}.
$$
The condition that $M \lambda  \in \Sigma_{3, \bar{2}}$ means each $a_j = 0$ or $1$
which imposes extra constraints on the $d_j$'s. \\

{\bf Claim 1.} {\em 
Suppose that $(d_0, d_1, ..., d_{2lm +k-1})$ with $0 \le k < m$ of $\lambda \in \sC(1,M)$
are fixed. Then at least one of the following conditions holds:\\

(i) There is at most one admissible value for $d_{2lm+k}$ in $\lambda~ (\bmod~3^{2lm+k+1})$.\\

(ii) There are two admissible values for $d_{2lm+k}$  for $ \lambda ~(\bmod~3^{2lm+k+1})$
and for any fixed choices
of 
$(d_{2lm+k+1}, d_{2lm+k+2}, ..., d_{(2l+1)m + k-1})$ at most three of
the four possible values of $(d_{2lm+k}, d_{(2l+1)m+k})$ give admissible sequences
for $\lambda ~(\bmod~3^{(2l+1)m+k}).$}\\

To prove the claim,  suppose that condition (i) doesn't hold. We then
examine the digit $a_{(2l+1)m+k}$ using 
\begin{eqnarray}\label{514a}
M \lambda  & \equiv&   b_0 d_{(2l+1)m +k}3^{(2l+1)m+k} +
b_m d_{(2l+m)+k} 3^{(2l+ 1)m+k}  \nonumber \\
&&~~~~~ +
M(\sum_{j=0}^{2lm+k-1} d_j 3^j) + b_0 d_{2lm +k}3^{2lm+k} ~(\bmod~3^{(2l+1)m +k+1}).
\end{eqnarray}
Define the digits $r_j$ by 
$$
M(\sum_{j=0}^{2lm+k-1} d_j 3^j) = \sum_{j=0}^{\infty} r_j 3^j, ~~~r_j \in  \{0,1,2\}.
$$
We assert that    \eqn{514a} then gives  the congruence 
\beql{515}
a_{(2l+1)m+k} \equiv b_0 d_{(2l+1)m+k} + b_m d_{2lm+k} + r_{(2l+l)m+k}  ~(\bmod~3).
\eeq
That is, we assert there cannot be any extra "carry digit"
from lower order terms that affects the $(2l+1)m+k$-th  $3$-adic digit,
coming from the addition of $b_0 d_{2lm+k} 3^{2m+k}$ in \eqn{514a}.
Namely, the extra term $b_0 d_{2lm+k} 3^k$, where
$d_{2lm+k}= 0$ or $1$ contributes nothing if $d_{2lm+k}=0$, while if $d_{2lm+k} =1$
By our assumption that (i) doesn't hold, 
both values $d_{2lm+k}=0, 1$ occur for  admissible $\lambda (\bmod~3^{2lm+k})$ for these digits. 
Since $b_0=1$ and the $3$-adic digit of
$M\lambda$ in the $(2lm+k+1)$-st place is $0$ or $1$, this digit  must have been $0$
when $d_{2lm+k}=0$, and $1$ when  $d_{2lm+k}=1$, so there can be no "carry digit"
in the addition of  $b_0 d_{2lm+k} 3^k$, as asserted.

Now consider the pairs $(d_{2lm+k}, d_{(2l+1)m+k})$. Of the four values $(00), (01), (10), (11)$
that these may take, the quantities $b_0 d_{(2l+1)m+k}+ b_{m} d_{2lm+k}$ with $b_0=1$
and $b_m=1$ or $2$ will cover all residue classes $(\bmod~3).$ In particular, at least
one choice will result in $a_{(2l+1)m+k} \equiv 2 ~(\bmod~3)$ in \eqn{515}, and so
give a non-admissible set of digits $(\bmod~3^{(2l+1)m+k})$. This proves (ii), and the
claim. \\


{\bf Claim 2.} {\em  For $M$ having the ternary expansion {\em \eqn{510}} and a given
$r \ge 2m$   there are
are at most $3^{\frac{1}{2}r + 2m}$ admissible congruence classes
in $\sC(1, M)$  $(\bmod~ 3^{r})$. }\\

To prove the claim, we
group the $3$-adic digits in
pairs  $(d_{2jm+k}, d_{(2j+1)m +k)})$, $ 0 \le k < m$, for all pairs with $(2j+1)m +k \le r$.
 There are at most $2m-1$ unpaired digits. 
Claim $1$ establishes that, conditional on
the choice of all other allowed digits, there are at most three permitted choices
for the set of paired digits.  For each unpaired digit there are at most two choices for
its value. Since the number of paired digits is at most $\frac{1}{2}(r+1)$
the total number of admissible sequences $(\bmod ~3^{r+1})$ is 
at most $3^{\frac{1}{2}(r+1)} 2^{2m-1}$,
which implies Claim 2. \\

To conclude the proof, Claim 2 implies that we have a covering $\sI_r$ of $\sC(1, M)$
with a set of at most $3^{(\frac{1}{2}r + 2m}$ sets, each of diameter $3^{-(r+1)}.$
For each $\epsilon>0$ this covering satisfies 
$$
\sum_{I \in \sI_r} |I|^{\frac{1}{2}+ \epsilon}  \le  3^{(\frac{1}{2}r + 2m} (3^{-(r+1)})^{\frac{1}{2} + \epsilon}
\le  3^{-(r+1) \epsilon}.
$$
Letting $r \to \infty$, this bound implies
$\dim_{H}(\sC(1, M))\le \frac{1}{2} + \epsilon.$ Letting $\epsilon \to 0$ gives the result. $~~~\bsq$\\

We do not know whether the bound in Theorem~\ref{th15} is sharp.  
However it is possible  to show
 that $\sC(1,7)$ has $ \dim_{H} \sC(1,7)= \log_3 (\frac{1+\sqrt{5}}{2}) \approx 0.43$.


%
%
\paragraph{Proof of Theorem~\ref{th17}.}
We suppose  are given $N$ a positive integer with $N \in \Sigma_{3, \bar{2}}\int \ZZ$
 and $1 \le M_1< M_2 < \cdots < M_k$ 
with all $NM_i \in \Sigma_{3, \bar{2}}.$ Our object is to obtain an explicit  nonzero lower bound
on the Hausdorff dimension $\dim_{H} (\sC(M_1, M_2, \cdots, M_k))$. 
We set $n$ equal to the number of ternary digits in $NM_k$, so that
$n= \lceil  \log_3 NM_k \rceil.$ Now we consider the set
$$
\Sigma_C:= \{ \lambda=( \cdots d_2 d_1 d_0)_3: ~ \mbox{all~blocks} ~~ d_{(k+1)n-1} \cdots 
d_{kn+1}d_{kn} \in \{  0^n ,~(N)_3 \}~\}
\subset  \Sigma_{3, \overline{2}}.
$$
Since each $NM_j \in  \Sigma_{3, \bar{2}} $ is an integer with at most $n$ ternary digits, we have
$$
M_j \Sigma_C:= \{ \lambda=( \cdots d_2 d_1 d_0)_3: ~ \mbox{all~blocks} ~
 d_{(k+1)n-1} \cdots d_{kn+1}d_{kn} \in \{ 0^n ,~ (NM_j)_3 \}~\}
\subset \Sigma_{3, \bar{2}}.
$$
Thus $\Sigma_C \subset \sC(M_1, M_2, \cdots , M_k)$. By inspection $\Sigma_C$ is
a Cantor set which has Hausdorff dimension
$$
\dim_{H} \Sigma_C = \frac{ \log_3 (2)}{\log_3 (3^n)} = \frac{\log_3(2)}{\lceil \log_3(NM_k) \rceil},
$$
and the result follows. $~~~\bsq$.\\


%
%
%
%

\section{ Furstenberg Conjecture and Transversality of Semigroup Actions}


In 1970 Furstenberg \cite[p. 43]{Fu70}
formulated  the following conjecture which is in the  same direction as
Erd\H{o}s's question.

\paragraph{Conjecture ${\bf 2}^{'}$.}{\it (Furstenberg) Suppose $p$
and $q$ are 
not powers of the same integer. Then the expansions to  the base $B=pq$
of the powers $\{ (p^n)_{pq}: n \ge 1\}$ 
have the property that 
any given  finite pattern of consecutive base $B$ digits occurs in
$(p^n)_{pq}$   for
all sufficiently large  $n$.} \\

\noindent For example, for $p=2$ and $q=3$,  this conjecture
asserts that  any given pattern
of base $B=6$ digits will occur as consecutive
digits  in  the base $6$ expansion of $(2^n)_6$, for all sufficiently large $n$.
The restriction to products $B=pq$ of two (or more) multiplicatively independent elements was
motivated by results in Furstenberg's seminal work  \cite{Fu67}. 
There he showed that for any irrational number $\theta$ the set
$\{ p^m q_n \theta (mod ~1): m, n \ge 0\}$ is dense on the torus
$\RR/\ZZ$. However it is well known that there is an 
uncountable set of  irrational numbers $\theta$
for which $\{ p^m \theta: m \ge 0\}$ is not dense on the torus.

Conjecture E in the introduction proposes nevertheless
that Furstenberg's conjecture continues to  hold
 when the base $B=q$ is a prime
(in the special case $p=2$, $q=3$).
More generally one can ask whether Furstenberg's conjecture might be
valid more generally for base $B$ expansions 
for arbitrary  $B$ with $gcd(B, p)=1$.

A main object of Furstenberg \cite{Fu70} was to introduce  a 
notion of transversality of two
semigroups of transformations $\sS_1$ and $\sS_2$ acting on
a compact metric space $X$ with respect to a (suitable)
dimension function $dim(A)$ defined on 
all closed sets $A$. 
 
 \begin{defi}~\label{de51}
 {\em Two closed sets $A$ and $B$ in a 
compact metric space $X$ are 
 {\em transverse} (for a given dimension function) if
$$
dim(A \cap B) \le \max( dim(A)+ dim(B) - dim(X), 0).
$$
}
\end{defi}

\begin{defi}~\label{de52}
{\em Two semigroups $\sS_1$ and $\sS_2$ acting on a compact metric space $X$ 
are {\em transverse} (for a given dimension function)
if  any closed $\sS_1$-invariant set $A$ 
and any closed $\sS_2$-invariant set $B$
are themselves transverse, for that dimension function.
}
\end{defi}

He obtained as an immediate consequence of this
definition the following result 
concerning simultaneous invariant sets (\cite[p. 42]{Fu70}),
which draws on earlier work (\cite{Fu67}).

\begin{prop}~\label{pr51}
(Furstenberg) Suppose that  $\sS_1$
and $\sS_2$ are transverse semigroups acting
on a compact metric space $X$, and  that  $\sS_1$ has
the additional property: \\

(*) If $A$ is a closed $\sS_1$-invariant set with
$dim(A) = dim(X)$, then $A=X$. \\

\noindent Then any proper closed subset of $X$ invariant
under both $\sS_1$ and $\sS_2$ has $dim(A)=0$.
\end{prop}

Furstenberg does not 
construct any transverse semigroups, but as evidence for their
existence shows  for the following pair of tranformation semigroups that their (nontrivial)
simultaneously invariant closed sets satisfy this property (\cite[Theorem 3]{Fu70}).

\begin{prop}~\label{pr52}
(Furstenberg) Let $\ZZ_r$ be the ring of $r$-adic integers, and suppose that
$r=pq$ with $p>1$ and $q>1$ not both powers of the same integer. Define
transformations $D_s(x)= \lfloor \frac{x}{s} \rfloor$, for $s= p, q,$ and $pq$, and note
that $D_{pq}= D_p D_q = D_q D_p$. Let $\sS_p$ and $\sS_q$ denote
the semigroups generated by $D_p$ and $D_q$, respectively.  If $A$ is a
simultaneously $\sS_p$ and $\sS_q$ invariant proper closed subset of $\ZZ_r$,
then  $A$  has Hausdorff dimension zero.
\end{prop}

\noindent The proof of this result draws on his earlier work (\cite{Fu67}).
Furstenberg \cite[p. 45] {Fu70} goes on to 
conjecture that $\sS_p$ and $\sS_q$ are transverse semigroups acting on $\ZZ_r$.

 Conjectures A and B in the introduction are partially motivated by  Furstenberg's
 framework but  fall outside it. One
could approach Conjecture A by  considering only the 
ternary expansions of fractional parts  $\{\{\lambda 2^n\}\}$, and thus
iterating $x \to 2x$ on the compact space $X= \RR/\ZZ$. This defines a  larger exceptional set
$\sE(\RR/\ZZ)$, which contains
$\sE(\RR)$.   Does $\sE(\RR/\ZZ)$ 
 have Hausdorff dimension zero? This set includes 
 all dyadic rationals (thus $\lambda=1$), which is a dense set
 in $\RR/\ZZ$, so its closure is the whole space $X$, and is not covered by
 Furstenberg's results.


Furstenberg's formulation does not apply to 
semigroups of transformations  on the real numbers
because $\RR$ is not compact. One may ask: {\em Can Furstenberg's framework be 
generalized to apply to semigroups
of operators acting on the real numbers, or the integers?}

%
%
%
%
\section{Concluding Remarks}


We conclude by reviewing some history  related to Erd\H{o}s's question. 
 Erd\H{o}s  \cite{Er79}  raised his question
on ternary expansions of $2^n$
in connection with his conjecture that the binomial coefficient
 ${\binom{2n}{n}}$ is not squarefree for all $n \ge 5$.
This binomial coefficient is divisible by $4$ 
except for  $n = 2^k$,
so it is natural to examine when larger primes divide
${\binom{2^{k+1}}{2^k}}$. Here one has 
$$
3~ \mbox{does~not~divide}~  {\binom{2^{k+1}}{2^k}} \Longleftrightarrow
 \mbox{The~ternary~ expansion ~of}~ 2^n~ \mbox{omits the digit}~2,
$$
as follows from Lucas's theorem (Lucas\cite{Lu78}, 
see Graham et al. \cite[Exercise 5.61]{GKP94}). This led 
Erd\H{o}s to raise his ternary expansion question,
since   a positive answer to it 
would establish his binomial coefficient conjecture. 

As it turned out,   Erd\H{o}s's  binomial coefficient conjecture  was
later  resolved  affirmatively,
 without answering the  ternary expansion question. In 1985  
Sarkozy \cite{Sa85} proved that ${\binom{2n}{n}}$
is not squarefree 
for all sufficiently large $n$. About 1995, 
Granville and Ramar\'{e} \cite{GR96}
and, independently,  Velammal \cite{Ve95} proved it for all $n \ge 5$.

The theme of this paper  is that Erd\H{o}s's unconventional question retains interest for
its own sake, even though the problem that originally motivated its study
has been solved. 

%
%
%
%


\noindent Jeffrey C. Lagarias \\
Dept. of Mathematics\\
The University of Michigan \\
Ann Arbor, MI 48109-1043\\
\noindent email: {\tt lagarias@umich.edu}
\end{document}